\documentclass[leqno]{amsart}
\usepackage{euscript,amssymb}
\usepackage[arrow,matrix]{xy}

\newtheorem{stuff}{Stuff}[section]
\newtheorem{theorem}[stuff]{\sl Theorem}
\newtheorem{proposition}[stuff]{\sl Proposition}
\newtheorem{lemma}[stuff]{\sl Lemma}
\newtheorem{corollary}[stuff]{\sl Corollary}

\newenvironment{definition}{%
\vskip1ex\refstepcounter{stuff}\trivlist \itemindent 0pt
\item[\hskip\labelsep\sl Definition \thestuff.]%
\ignorespaces}{\endtrivlist\vskip1ex}%
\newenvironment{remark}{%
\vskip1ex\refstepcounter{stuff}\trivlist \itemindent 0pt
\item[\hskip\labelsep\sl Remark \thestuff.]%
\ignorespaces}{\endtrivlist\vskip1ex}%

\newtheorem{s-proposition}[sstuff]{\sl Proposition}
\newtheorem{s-lemma}[sstuff]{\sl Lemma}
\newtheorem{s-corollary}[sstuff]{\sl Corollary}
\newenvironment{s-definition}{%
\vskip1ex\refstepcounter{sstuff}\trivlist \itemindent 0pt
\item[\hskip\labelsep\sl Definition \thesstuff.]%
\ignorespaces}{\endtrivlist\vskip1ex}%
\newenvironment{s-remark}{%
\vskip1ex\refstepcounter{sstuff}\trivlist \itemindent 0pt
\item[\hskip\labelsep\sl Remark \thesstuff.]%
\ignorespaces}{\endtrivlist\vskip1ex}%

\let\rar\rightarrow
\let\lar\longrightarrow
\let\llar\longleftarrow
\let\dar\downarrow
\let\sear\searrow
\let\swar\swarrow
\let\near\nearrow

\let\hra\hookrightarrow
\let\mt\mapsto
\let\lmt\longmapsto
\let\xrar\xrightarrow

\font\tenmsa=msam10 %
\newcommand\hdashpiece{%
{\vrule height2.75pt depth-2.35pt width2.3pt \kern1.7pt}}%
\newcommand\hdashpieces{%
{\hdashpiece\hdashpiece\hdashpiece\hdashpiece}}%
\newcommand\dashto{\mathrel{%
\hdashpiece\hdashpiece\kern-0.4pt\hbox{\tenmsa K}}}%
\newcommand\dashar{\mathrel{%
\hdashpieces\kern-0.4pt\hbox{\tenmsa K}}}%
\newcommand\longdashar{\mathrel{%
\hdashpieces\hdashpieces\kern-0.4pt\hbox{\tenmsa K}}}%
\newcommand\ldashar{\mathrel{%
\hbox{\tenmsa L}\hdashpieces}}%
\newcommand\lrdash{\mathrel{%
\hbox{\tenmsa L}\kern1pt\hdashpiece\kern-0.4pt\hbox{\tenmsa K}}}
\newcommand\lrdashar{\mathrel{%
\hbox{\tenmsa L}\kern1pt\hdashpieces\kern-0.4pt\hbox{\tenmsa K}}}
\newcommand\longlrdashar{\mathrel{%
\hbox{\tenmsa L}\kern1pt\hdashpieces\hdashpieces
\kern-0.4pt\hbox{\tenmsa K}}}%
\font\arrowhead=cmex10
\def\downarrowhead{\hbox{\arrowhead\char '171}} 
\newcommand\vdashpiece{%
{\kern-0.2pt{\vrule height2.3pt depth0pt width0.4pt}}} %
\newcommand\vdashto{{\begin{array}{c}
    \vdashpiece\\[-8pt]\vdashpiece\\[-11pt]\downarrowhead
    \end{array}}}%
\def\vequal#1{{\begin{array}{cc} 
\vrule height7pt depth #1 width0.3pt\kern-1ex
&\kern-0.9ex\vrule height7pt depth #1 width0.3pt
\end{array}}} 

\let\euf\EuScript 
\let\cal\mathcal
\let\mbb\mathbb

\let\mfrak\mathfrak
\DeclareFontFamily{OT1}{rsfs}{}
\DeclareFontShape{OT1}{rsfs}{n}{it}{<->rsfs10}{}
\DeclareMathAlphabet{\curly}{OT1}{rsfs}{n}{it}
\let\crl\curly

\let\ovl\overline
\let\unl\underline
\let\tld\tilde
\let\wtld\widetilde
\let\wht\widehat
\let\nit\noindent
\let\disp\displaystyle
\let\srel\stackrel
\let\vphi\varphi
\let\veps\varepsilon
\newcommand\dbar{{\bar{\partial }}}
\newcommand\ort{\mathrel{{\vrule width4.0pt height0.4pt depth0pt
                  \vrule width0.4pt height6.0pt depth0pt\,}}}
\newcommand\rd{{\rm d}}
\newcommand\Hom{\mathop{\sf Hom}\nolimits}

\newcommand\invq{{\slash\kern-0.65ex\slash}}
\newcommand\EG{{\rm EG}}
\newcommand\BG{{\rm BG}}
\newcommand\EK{{\rm EK}}
\newcommand\BK{{\rm BK}}
\newcommand\sev{{\hbox{\sc ev}}}
\newcommand\rmr{{\rm r}}
\newcommand\rmR{{\rm R}}
\newcommand\s{{\rm s}}
\newcommand\sst{{\rm ss}}
\numberwithin{equation}{section}
\begin{document}
\title{GW invariants and invariant quotients}
\author{Mihai Halic}
\address{Fakult\"at f\"ur mathematik NA\\ 
Ruhr-Universit\"at Bochum\\ 
44780 Bochum\\ Germany}
\email{halic@ag.ruhr-uni-bochum.de}
\maketitle
\markboth{GW INVARIANTS AND INVARIANT QUOTIENTS}{MIHAI HALIC}

\section*{Introduction}

Gromov-Witten theory became a major tool in enumerative geometry 
because the Gromov-Witten invariants (GW-invariants for short) 
give in some cases the number of curves in projective varieties, 
satisfying certain incidence conditions. However, it is usually 
quite difficult to explicitely compute these invariants, and 
therefore it is useful to know their behavior under `modifications' 
of the variety we start with; typical examples are the blow-up or 
the symplectic connected sum.

In the present paper we are dealing with group actions. The starting 
point of this study was the question: given a projective variety and 
a group acting on it, is there any relationship between the 
GW-invariants of this variety and those of its quotient for the group 
action? In this form, it is rather hopeless to answer the question, so 
that it eventually became: given a projective variety $X$ with a very 
ample line bundle ${\cal O}_X(1)\rar X$ and a connected, linearly 
reductive group $G$ whose action on $X$ is linearized in 
${\cal O}_X(1)$, is there any relationship between the GW-invariants 
of $X$ and those of its invariant quotient $X\invq G$, under the 
additional assumption that $G$ acts {\em freely} on the 
$G$-semi-stable locus of $X$? Certainly, such an assumption makes the 
geometric invariant theory on $X$ trivial because in this case 
$X^\sst\rar X\invq G$ is simply a principal $G$-bundle. We shall see 
however, that this assumption appears quite naturally in the context, 
and even so we will have to face rather complicated situations.

The strategy adopted to attack the question above is the following: 
the $G$-action on $X$ induces one on the space of stable maps to $X$ 
and two maps in the same $G$-orbit, with image contained in the 
semi-stable locus of $X$, induce the same map into the quotient. 
The first idea which comes to mind is to compare the invariant 
quotient, for this action, of the moduli space of stable maps to $X$ 
with the moduli space of stable maps to the quotient $X\invq G$. 
For a reason which will become clear in a moment, this strategy is 
correct only in genus zero, for in higher genera changes are needed. 

Let me present more precisely the setting of the problem: the group 
$G$ and the variety $X$ are as above, and $A\in H_2(X;{\mbb Z})$ is a 
homology class which can be represented by a morphism $C\rar X^\sst$, 
with $C$ a smooth projective curve of genus $g$; we let 
$\hat A:=\phi_*A\in H_2(X\invq G;{\mbb Z})$ be the push-forward 
for the projection $\phi:X^\sst\rar X\invq G$. What we would like 
to compute are the GW-invariants of $X\invq G$ corresponding to 
the class $\hat A$.

For making clear the logic of the article, I shall start with a 
naive comparison of the expected dimensions of the moduli spaces 
of stable maps involved here (the precise meaning of the notations 
will be given later on): 
$$
\begin{array}{l}
\phantom{x}
\hat D:=\text{exp.dim\,}\ovl M_{g,k}(X\invq G,\hat A)
=(3-\dim X\invq G)(g-1)+c_1(X\invq G)\cdot\hat A+k,\\[1,5ex] 
\begin{array}{ll}
D-\dim G\kern-1ex&:=\exp.\dim.\ovl M_{g,k}(X,A)\invq G\\ 
&\;=(3-\dim X)(g-1)+c_1(X)\cdot A+k-\dim G,
\end{array}
\end{array}
$$
and therefore

$(\star)$
\phantom{x}\hfill 
$\hat D - (D - \dim G) = g\cdot\dim G.$
\hfill\phantom{x}

\nit From this computation we deduce that in general the space 
$\ovl M_{g,k}(X\invq G,\hat A)$ is {\em larger} than 
$\ovl M_{g,k}(X,A)\invq G$, the only exception happening in genus 
zero. Very shortly, the explanation for this phenomenon is that 
the projective line is the only one smooth curve which has the 
property that a topologically trivial, holomorphic principal 
$G$-bundle over it is also holomorphically trivial. In higher 
genera, holomorphic principal $G$-bundles with fixed topological 
type depend on `moduli', whose number agrees with the difference 
$(\star)$ above. This is the reason why for computing higher genus 
invariants of $\hat X$ we will need to consider maps into a 
larger variety $\bar X$ whose construction, in the case when 
$G$ is a torus, is given in lemma \ref{lm:u-bar}.

The article is organized as follows: the first section recalls some 
basic facts about stable maps and their moduli spaces, the reference 
being \cite{FP}. 

In section \ref{sct:ss-points} we describe the $G$-semi-stable 
points of the moduli space of stable maps $\ovl M_{g,k}(X,A)$. 
The results obtained in this section hold in full generality, no 
matter what the $G$-action on $X$ looks like. We obtain the sufficient 
result (theorem \ref{thm:ss-maps1}) which says that a map with image 
contained in the semi-stable locus of $X$ is  $G$-semi-stable as a 
point of $\ovl M_{g,k}(X,A)$ and a necessary result 
(corollary \ref{cor:case-torus}) which says that a stable map 
representing a $G$-semi-stable point of $\ovl M_{g,k}(X,A)$ does 
not have its image contained in the unstable locus of $X$. 

Section \ref{sct:symp-pers} characterizes the semi-stable points of 
$\ovl M_{g,k}(X,A)$ from a symplectic point of view, which will be 
useful later on in section \ref{sct:spaces} where we will give an 
algebro-geometric construction of the space of maps needed for defining 
certain `Hamiltonian invariants'. We compute an explicit formula for 
the moment map on the space of stable maps which corresponds to a 
${\mbb C}^*$-action (proposition \ref{prop:mm-on-global-cover}), and 
using it we give (theorem \ref{thm:geom-picture}) a second proof for 
theorem \ref{thm:ss-maps1}.

Finally, section \ref{sct:conseq1} closes the first part of the 
article giving a partial answer to the initial problem, that of 
comparing the genus zero GW-invariants of $X$ and $X\invq G$. 
Theorem \ref{thm:comparison} states, under certain transversality 
assumptions which are technical in nature, that if $G$ acts freely 
on the semi-stable locus of $X$, $A$ is a spherical class for which 
a representative may be found to lie entirely in the stable locus 
$X^\sst$ and $\hat A$ denotes the push-forward class in $X\invq G$, 
then 
$$
GW_{X\invq G,\hat A}^{0,k}(\hat a_1,\dots,\hat a_k) 
=GW_{X,A}^{0,k}(a_1\zeta,\dots,a_k).
$$
In this equality the $a_i$'s are the classes on $X$ determined by 
the $\hat a_i$'s on $X\invq G$ {\it via} the rational quotient map 
(essentially by pull-back), and $\zeta$ denotes the class of a 
rational slice for $X\dashto X\invq G$. We conclude the section 
with some explicit computations. 

Sections \ref{sct:conseq2}, \ref{sct:spaces} and \ref{sct:the-inv} 
grew out from the attempt at understanding from an algebraic point 
of view the Hamiltonian invariants defined in \cite{cgs} and 
\cite{mun}, in the case of torus actions. Maybe I should say a 
word about the origin of this interest: on one hand, in sections 
\ref{sct:ss-points} and \ref{sct:symp-pers} we have tried to develop 
the algebraic tools needed for computing the GW-invariants of a 
quotient variety while, on the other hand, the papers \cite{cgs} 
and \cite{mun} bring into the scene some new invariants, constructed 
using real-analytic methods, which are associated to Hamiltonian group 
actions on symplectic manifolds. For particular choices of certain 
parameters, these ones should compute the GW-invariants of the 
Marsden-Weinstein quotient of the symplectic variety on which the 
group acts. Since in the case of projective varieties the symplectic 
reduction and the geometric invariant quotient basically agree, it was 
natural to try relating the two constructions: the algebro-geometric 
and the real-analytic one. 

It turns out that the construction of these `Hamiltonian invariants' 
fits well in the frame of the previous sections. We identify 
(proposition \ref{prop:1-1}) the moduli spaces introduced in 
\cite{cgs,mun} with the space of stable maps into the variety 
$\bar X:={\cal P}\times_{{({\mbb C}^*)}^r}X$, where 
${\cal P}\rar{(Pic^0C)}^r\times C$ denotes a Poincar\'e bundle 
parameterizing principal $G={({\mbb C}^*)}^r$-bundles over the 
curve $C$, so that the real-analytic and the algebraic points of 
view actually coincide.

In the last section we prove (theorem \ref{thm:compar2}), under the 
same transversality assumptions as in theorem \ref{thm:comparison}, 
the conjecture formulated in \cite{cgs} which states that for 
certain choices of the parameters these Hamiltonian invariants of 
$X$ coincide with Gromov-Witten invariants on the quotient 
$X\invq G$.

\section{Some preparatory material}
{\label{sct:reminders}}

In this section I shall recall from \cite{FP} the definition of 
a stable map, the construction of the space $\ovl M_{g,k}(X ,A)$ 
as a projective scheme and the description of an ample line bundle 
on it. 

In the whole paper $X$ denotes a complex projective, irreducible 
variety and ${\cal O}_X(1)\rar X$ denotes a very ample line bundle 
on it; let $X\hra {\mbb P}^r$ be the embedding given by the linear 
system of ${\cal O}_X(1)$.  

\begin{definition}{\label{def:stable-map}} 
A {\em stable map} $[(C ,\unl x,u)]$ to $X$ consists of the 
following data: 

(i) a connected, reduced, complete algebraic curve $(C,\unl x)$ 
with $k$ distinct marked points $\unl x=(x_1,\dots ,x_k)$. The 
singularities of $C$ are at most ordinary double points and the 
markings lie in the smooth locus of $C$; 

(ii) the equivalence class of a morphism $u:(C,\unl x)\rar X$. 
Two morphism $u:(C,\unl x)\rar X$ and $u':(C',\unl x')\rar X$ 
are equivalent if there exists an isomorphism $\gamma :C\rar C'$ 
such that $\gamma(x_j)=x_j'$ for $j=1,\dots,k$ and 
$u=u'\circ\gamma$; 

(iii) the stability means that the automorphism group of the map 
$(C,\unl x,u)$ is finite. 
\end{definition}

By abuse of language, we shall often call a stable map and denote 
it by $(C,\unl x,u)$, the morphism $u:(C,\unl x)\rar X $ itself 
(satisfying (i)+(iii)) and not the equivalence class defined by it.  

A stable map $[(C,\unl x,u)]$ is said to represent the 
$2$-homology class $A\in H_2(X;{\mbb Z})$ if $u_*[C ]=A$. 
Composing such a map with the inclusion $X\hra {\mbb P}^r$, one 
obtains a stable map to ${\mbb P}^r$ which represents an integral 
multiple $d$ of the class $\ell$ of a line in ${\mbb P}^r$. In the 
sequel we shall briefly recall the construction of the moduli space 
$\ovl M_{g,k}({\mbb P}^r,d)$ of stable maps of genus $g$, with $k$ 
marked points and representing the class 
$d\ell\in H_2({\mbb P}^r;{\mbb Z})$. The moduli space 
$\ovl M_{g,k}(X ,A)$ will be a closed subscheme of it. 

It is immediate that the stability condition is equivalent 
to the fact that the line bundle $L_{(C,\unl x,u)}:=
\omega_C(x_1+ \dots +x_k)\otimes u^*{\cal O}_X(3)\rar C$ is ample, 
and a combinatorial argument proves that there is an integer 
$f\geq 1$ having the property that for any stable map $(C,\unl x,u)$ 
to ${\mbb P}^r$ of genus $g$ and with $k$ marked points, 
representing the class $d\ell\in H_2({\mbb P}^r;{\mbb Z})$, the 
line bundle $L_{(C ,\unl x,u)}^{\otimes f}\rar C$ is very ample and 
$H^1(C,L_{(C,\unl x,u)}^{\otimes f})=0$. Fix such an integer $f$ 
once for all. For any $k$-pointed, genus $g$ stable map 
$(C,\unl x,u)$ representing the class $d\ell$ one gets an embedding 
$$
u_Y:=(|L_{(C ,\unl x,u)}^{\otimes f}|,u):C\lar  
{\mbb P}\left( {H^0(C,L_{(C ,\unl x,u)}^{\otimes f})}^\vee \right)
\times {\mbb P}^r\cong {\mbb P}^N\times {\mbb P}^r=:Y,
$$ 
into a product of two projective spaces. 
The dimension of the projective spaces 
${\mbb P}\left( {H^0(C ,L_{(C ,\unl x,u)}^{\otimes f})}^\vee \right)$ 
is independent of the choice of the stable map with the properties 
mentioned above and is given by a Riemann-Roch formula. Notice that 
$C$ is determined as a subvariety of $Y$ up to a $PGl(N+1)$-action 
and that the Hilbert polynomial $P$ of $C$ inside $Y$ does not depend 
on the choice of the stable map $(C,\unl x,u)$. Denoting by 
${\cal H}:={\cal H}ilb_Y^P$ the Hilbert scheme of closed subschemes 
of $Y$ whose Hilbert polynomial is $P$, to each stable map 
$(C,\unl x,u)$ one associates a point in ${\cal H}\times Y^k$ as 
follows:
$$
(C ,\unl x,u)\mt (u_{Y*}C ,u_Y(x_1),\dots ,u_Y(x_k)).
$$ 
The natural $PGl(N+1)$-action on ${\mbb P}^N$ induces an action on 
${\cal H}\times Y^k$ and two stable maps $(C_{1,2},\unl x_{1,2}, 
u_{1,2})$ are isomorphic if and only if they are in the same 
$PGl(N+1)$-orbit. The stability condition translates into the fact 
that the stabilizer of any stable map under this action is finite. 

It is proved in \cite{FP} that there is a certain subscheme ${\cal S}$ 
of ${\cal H}\times Y^k$ such that $\ovl M_{g,k}({\mbb P}^r,d)
={\cal S}/PGl(N+1)$. One of the main results obtained in that paper is 
that $\ovl M_{g,k}({\mbb P}^r,d)$ is a separated and proper scheme, 
projective over $\mbb C$. Following \cite{FP}, we are going to describe 
an ample line bundle on $\ovl M_{g,k}({\mbb P}^r,d)$. The Hilbert 
scheme $\cal H $ is projective, an embedding of it into a projective 
space being obtained as follows: ${\cal L}:={\cal O}_{{\mbb P}^N}(1)
\boxtimes {\cal O}_{{\mbb P}^r}(1)\rar Y$ is a very ample line 
bundle and for large enough integral values of $l$, the restriction 
homomorphism
$$
\kern-0.5pt{\euf W} :=H^0(Y,{\cal L}^l)^{(0)}\oplus 
\bigoplus _{j=1}^kH^0(Y,{\cal L} ^l)^{(j)}
\lar H^0(Y,{\cal O}_C \otimes {\cal L} ^l)\oplus 
\bigoplus _{j=1}^k {\cal L}_{x_j}^l
$$ 
is surjective at each point $(C,x_1,\dots ,x_k)\in {\cal H}\times 
Y^k$. For distinguishing between the different direct summands of 
${\euf W}$, we have used upper indices for the same vector space 
$H^0(Y,{\cal L}^l)$. Moreover, under the same assumption that $l$ 
is large enough, the dimension 
$$
{\rm dim}\,H^0(Y,{\cal O}_C \otimes {\cal L}^l)\oplus 
\bigoplus _{j=1}^k {\cal L}_{x_j}^l=P(l)+k=:q+k
$$ 
is independent of the point $(C,x_1,\dots,x_k)\in{\cal H}\times Y^k$. 

It is proved in \cite{Gk} that the map 
$$
\begin{array}{ccl} 
{\cal H}\times Y^k & 
\lar & 
Gr_{q+k}({\euf W}^\vee ),\\[1.5ex] 
(C,\unl x) & 
\mt & \disp 
{\left( H^0(Y,{\cal O}_C \otimes {\cal L}^l)\oplus 
\bigoplus _{j=1}^k {\cal L}_{x_j}^l \right) }^\vee 
\end{array}
$$
is a closed immersion. Composing it with the usual projective 
embedding of the Grassmann variety, one obtains the projective 
embedding of ${\cal H}\times Y^k$: 
$$
\begin{array}{ccl} 
{\cal H}\times Y^k & 
\lar & 
\disp 
{\mbb P}\left( \bigwedge ^{q+k}{\euf W}^\vee\right),\\[1.5ex]  
(C,\unl x) & 
\mt & 
\disp 
{\rm det}\,\left({\left( H^0(Y,{\cal O}_C \otimes {\cal L}^l)
\oplus \bigoplus _{j=1}^k {\cal L}_{x_j}^l \right) }^\vee \right).
\end{array}
$$
The induced ample line bundle on ${\cal H}\times Y^k$ is 
\begin{align}{\label{ample-line-bdl}}
{\rm det}\,{\cal Q}_k :={\rm det}\,\left( {\cal Q}
\boxplus{({\cal L}^l)}^{\boxplus k}\right),
\end{align}
where ${\cal Q}\rar Gr_q\left( {H^0(Y,{\cal L}^l)}^\vee \right)$ is 
the universal quotient bundle. The fibre of ${\cal Q}_k$ at a point 
$(C,\unl x)\in {\cal H}\times Y^k$ is 
${{\cal Q}_k}_{(C,\unl x)}=H^0(Y,{\cal O}_C \otimes {\cal L}^l)\oplus 
\bigoplus _{j=1}^l{\cal L}_{x_j}^l$. Since ${\cal Q}_k$ is invariant 
under the $PGl(N+1)$-action, det$\,{\cal Q}_k$ is also. It is proved 
in \cite{FP} that det$\,{\cal Q}_k $ descends to an {\em ample} line 
bundle on $\ovl M_{g,k}(X,A)={\cal S}/PGl(N+1)$. 


\section{The semi-stable points on $\ovl M_{g,k}(X,A)$} 
{\label{sct:ss-points}}

Here is the setup: $G$ is a connected, linearly reductive, complex 
algebraic group which acts on the complex, irreducible, projective 
variety $X$. The action is linearized in the very ample line bundle 
${\cal O}_X(1)\rar X $. This action naturally induces one on the 
space of stable maps by
$$
G\times \ovl M_{g,k}(X ,A)\lar \ovl M_{g,k}(X ,A),
$$ 
$$
\left( g,[(C ,\unl x,u)]\right)\mt [(C,\unl x,gu)],
$$  
where $gu:C \rar X $ is defined by $(gu)(p):=g\cdot u(p)$ for all 
$p\in C $. In order to compute the geometric invariant quotient of 
$\ovl M_{g,k}(X,A)$, we need a linearization of the action in an 
ample line bundle.

The linearized $G$-action on $X$ extends to a linearized action on 
$$
{\cal O}_{{\mbb P}^r}(1)\rar
{\mbb P}^r:={\mbb P}\left( {H^0(X,{\cal O}_X(1))}^\vee\right)
$$ 
such that $X$ is invariant. We have already mentioned in the 
previous section that $\ovl M_{g,k}(X,A)$ is a closed subscheme 
of $\ovl M_{g,k}({\mbb P}^r,d)$ and on this last one we have 
described an ample line bundle. Our next task is to linearize 
the $G$-action in it. 
 
Using the notations of the previous section, 
$\ovl M_{g,k}({\mbb P}^r,d)={\cal S}/PGl(N+1)$, where $\cal S$ 
is some subscheme of ${\cal H}\times Y^k$. On ${\cal H}\times Y^k$ 
there are two actions: the first one is the $PGl(N+1)$-action on 
${\cal H}$ and the second one is the $G$-action induced by that on 
the factor ${\mbb P}^r$ in $Y={\mbb P}^N\times {\mbb P}^r$. Since 
the $PGl(N+1)$ and the $G$-actions on $Y$ commute, the induced 
actions on ${\cal H}\times Y^k$ commute also. For this reason, 
$G$-semi-stable points on $\ovl M_{g,k}({\mbb P}^r,d)$ will be the 
images of $G$-semi-stable points of $\cal S$; therefore it is 
enough to describe the linearized $G$-action on 
det$\,{\cal Q}_k\rar {\cal H}\times Y^k$. Since $G$ acts on 
${\cal L}\rar Y$, it acts also on $H^0(Y,{\cal L}^l)$ by 
${(gS)}_y:=gS_{g^{-1}y}$ for all $S\in H^0(Y,{\cal L}^l)$ and all 
$y\in Y$. The dual action on ${H^0(Y,{\cal L}^l)}^\vee $ is given by 
$(g,\Sigma )\mapsto g\Sigma$, where 
$\langle g\Sigma ,S\rangle:=\langle \Sigma ,g^{-1}S\rangle $. 
The induced action on $\disp\bigwedge^{q+k}{\euf W}^\vee$ 
is now  obvious because ${\euf W}= H^0(Y,{\cal L}^l)^{(0)}\oplus 
\bigoplus _{j=1}^kH^0(Y,{\cal L}^l)^{(j)}$. 

For finding the $G$-semi-stable points we shall use the
Hilbert-Mumford criterion. Given $\lambda :{\mbb C}^*\rar G$ a 
one parameter subgroup of $G$ (a 1-PS for short), there is a (finite) 
direct sum decomposition corresponding to the characters of 
${\mbb C}^*$: 
$$
H^0(Y,{\cal L}^l)=\bigoplus _{m\in {\mbb Z}}H^0(Y,{\cal L}^l)_m,
$$ 
$$
\lambda (t)S=t^mS \kern1em \forall t\in {\mbb C}^*\kern1em 
\forall S\in H^0(Y,{\cal L}^l)_m.
$$ 
We want to find out when the point 
$$
\left[ 
\bigwedge ^{q+k}{ \left( H^0(Y,{\cal O}_C\otimes {\cal L}^l)^{(0)}
\oplus \bigoplus _{j=1}^k{\cal L}_{x_j}^l \right) }^\vee \right]\in 
{\mbb P}\left( \bigwedge ^{q+k}{\euf W\;}^\vee \right)
$$ 
is $\lambda $-semi-stable, so we have to study the ${\mbb C}^*$-orbit 
of a representative of this point in 
$\disp\bigwedge^{q+k}{\euf W}^\vee$. 

Let $\sigma_1, \dots ,\sigma_q$ be a basis of 
$\Hom_{\mbb C}\left( H^0(C,{\cal L}^l)^{(0)},{\mbb C}\right)$ and 
$\tau_1,\dots,\tau_k$ be generators of 
$\Hom({\cal L}_{x_j}^l,{\mbb C})$, $j=1,\dots,k$. 
Notice that the choice of the $\tau_j$'s is equivalent to the choice 
of representatives $x_j'\in {\mbb C}^N\times{\mbb C}^{r+1}$ of 
$x_j=(x_{j,1},x_{j,2})\in Y$ because ${\cal L}_{x_j}=
{\cal O}_{{\mbb P}^{N-1}}(1)_{x_{j,1}}\otimes 
{\cal O}_{{\mbb P}^r}(1)_{x_{j,2}}$. 
Using the epimorphism 
$$
{\euf W}=H^0(Y,{\cal L}^l)^{(0)}\oplus \bigoplus_{j=1}^k 
H^0(Y,{\cal L}^l)^{(j)}\stackrel{\mbox{$\imath _{_C }$}}{\lar}
H^0(C ,{\cal L}^l)\oplus\bigoplus_{j=1}^k {\cal L}_{x_j}^l\lar 0,
$$ 
$\sigma_1,\dots,\sigma_q,\tau_1,\dots,\tau_k$ can be extended to 
linear functionals on ${\euf W}$ 
$$
{\euf S}_1, \dots ,{\euf S}_q :H^0(Y,{\cal L}^l)^{(0)}
\longrightarrow {\mbb C}
$$ 
$$
\langle {\euf S} _j,S\rangle :=\langle \sigma_j, 
\imath _{_C }S\rangle \kern1em j=1,\dots ,q \kern1em  
\forall S\in H^0(Y,{\cal L}^l)^{(0)} 
$$ 
and 
$$ 
{\euf T}_1, \dots ,{\euf T}_k:H^0(Y,{\cal L}^l)\lar {\mbb C}
$$ 
$$
\langle {\euf T} _j,S\rangle :=\langle \tau _j,S(x_j)\rangle 
\kern1em j=1,\dots ,k \kern1em \forall S\in  H^0(Y,{\cal L}^l)^{(j)}.
$$ 
The ${\euf T}_j$'s represent just the evaluations of the 
homogeneous polynomial $S$ at the points $x_j'$ representing $x_j$. 

Let us remark that the linear functionals ${\euf S}_j$ act only on 
$H^0(Y,{\cal L}^l)^{(0)}$ and evaluate identically to zero on the 
other copies  $H^0(Y,{\cal L}^l)^{(j)}$, $j\neq 0$. A similar remark 
is valid for the ${\euf T} _j$'s: they evaluate identically to zero 
on $H^0(Y,{\cal L}^l)^{(j')}$, $j'\neq j$. 

The semi-stability condition reads 
$$
0\not \in {\ovl{{\mbb C}^* \cdot {\euf S}_1\wedge\dots 
\wedge{\euf S} _q\wedge 
{\euf T} _1\wedge \dots \wedge {\euf T} _k }}^
{\bigwedge ^{q+k}{\euf W\;}^\vee }
$$ 
which is equivalent to the existence of $S_1,\dots,S_{q+k},S'_1,\dots,
S'_{q+k}\in {\euf W}$  such that 
\begin{align}{\label{ss-condition}}
\left\{
\begin{array}{l}
\disp 
0\neq\lim_{t\rar 0}\langle
\lambda(t) \kern-2pt\cdot\kern-2pt \left( 
{\euf S}_1\wedge \dots \wedge {\euf S} _q\wedge 
{\euf T}_1 \wedge \dots \wedge {\euf T} _k \right), 
S_1\wedge \dots \wedge S_{q+k}\rangle,\\[1.5ex]
\disp 
0\neq\lim_{t\rar\infty}\langle
\lambda(t) \kern-2pt\cdot\kern-2pt \left( 
{\euf S}_1\wedge \dots \wedge {\euf S}_q\wedge 
{\euf T}_1 \wedge \dots \wedge {\euf T}_k \right), 
S'_1\wedge \dots \wedge S'_{q+k}\rangle.
\end{array}
\right.
\end{align} 
Each vector $S_j,S'_j$ is the sum of $1+k$ vectors corresponding to
the direct sum decomposition of $\euf W $. Moreover, as pointed out 
before, each of the ${\euf S}_j\,$'s and ${\euf T}_j\,$'s evaluate 
non-zero only on vectors in a certain component of $\euf W\,$. 
Consequently, for this last condition to be satisfied, one may assume 
that 
$$
S_1,S'_1\dots, S_q,S'_q\in H^0(Y,{\cal L}^l)^{(0)} 
\kern1em {\rm and}\kern1em 
S_{q+j},S'_{q+j}\in H^0(Y,{\cal L}^l)^{(j)}\ \ j=1,\dots ,k.
$$ 

Since the ${\mbb C}^*$-action of the 1-PS $\lambda $ of $G$ induces 
the decomposition $H^0(Y, {\cal L}^l)\break =\oplus _{m\in {\mbb Z}}\, 
H^0(Y,{\cal L}^l)_m$, we can further assume that 
$$
S_j\in H^0(Y,{\cal L}^l)_{m_j}^{(0)},\quad 
S'_j\in H^0(Y,{\cal L}^l)_{m'_j}^{(0)}
\kern1em {\rm for }\kern1em  
j=1,\dots ,q
$$ 
and 
$$
S_{q+j}\in H^0(Y,{\cal L}^l)_{m_{q+j}}^{(j)},\quad 
S'_{q+j}\in H^0(Y,{\cal L}^l)_{m'_{q+j}}^{(j)}
\kern1em {\rm for}\kern1em 
j=1,\dots ,k.
$$ 
We are now going to compute the first condition in \eqref{ss-condition}:
$$
\begin{array}{l}  
\langle 
\lambda (t)\kern-2pt\cdot\kern-2pt\left( 
{\euf S} _1\wedge \dots \wedge {\euf S} _q\wedge 
{\euf T} _1\wedge \dots \wedge {\euf T} _k \right)\kern-1pt , 
S_1\wedge \dots S_q\wedge S_{q+1}\wedge \dots \wedge S_{q+k}
\rangle \\[1em] 

= \langle 
{\euf S} _1\wedge \dots \wedge {\euf S} _q\wedge  
{\euf T} _1\wedge \dots \wedge {\euf T} _k , \\[1ex]

\kern1.5cm  
\lambda (t^{-1})S_1\wedge \dots \wedge \lambda (t^{-1})S_q \wedge 
\lambda (t^{-1})S_{q+1}\wedge \dots \wedge \lambda (t^{-1})S_{q+k} 
\rangle \\[1em]  

=\langle 
{\euf S} _1\wedge \dots \wedge {\euf S} _q\wedge  
{\euf T} _1\wedge \dots \wedge {\euf T} _k ,\\[1ex] 

\kern1.5cm  
t^{-m_1}S_1\wedge \dots \wedge t^{-m_q}S_q\wedge 
t^{-m_{q+1}}S_{q+1}\wedge \dots \wedge t^{-m_{q+k}}S_{q+k} 
\rangle \\[1em] 

=t^{-\sum _{j=1}^{q+k}m_j}\cdot 
\langle 
{\euf S} _1\wedge \dots \wedge {\euf S} _q,
S_1\wedge \dots \wedge S_q 
\rangle 
\cdot 
\prod _{j=1}^k\langle {\euf T} _j,S_{q+j}
\rangle \\[1em] 
=t^{-\sum _{j=1}^{q+k}m_j}\cdot 
\langle 
\sigma _1\wedge \dots \wedge \sigma _q, 
\imath _{_C }^*S_1\wedge \dots \wedge \imath _{_C }^*S_q
\rangle 
\cdot 
\prod _{j=1}^kS_{q+j}(x_j').
\end{array}
$$
Doing the computations corresponding to the second condition 
in \eqref{ss-condition} we find the 

\begin{proposition}{\label{prop:ss-maps}} 
The point $[(C ,\unl x,u)]\in \ovl M_{g,k}(X ,A)$ is $G$-semi-stable
if and only if for any {\rm 1-PS} $\lambda :{\mbb C}^*\rar G$ there 
are sections 
$$
S_1,S'_1\dots ,S_q,S'_q,
S_{q+1},S'_{q+1},\dots ,S_{q+k},S'_{q+k}
\in H^0(Y,{\cal L}^l)
$$ 
satisfying the properties: 

{\rm (i)}  
\begin{tabular}[t]{l}
$\lambda (t)S_j=t^{m_j}S_j$ for $j=1,\dots ,q+k$ with 
$\disp \sum _{j=1}^{q+k}m_j\geq 0;$\\ 
$\lambda (t)S'_j=t^{m'_j}S_j$ for $j=1,\dots ,q+k$ with 
$\disp \sum _{j=1}^{q+k}m'_j\leq 0;$
\end{tabular}
\smallskip

{\rm (ii)}  
$\{ \imath_{_C }^*S_1,\dots ,\imath _{_C}^*S_q\}$  and 
$\{ \imath_{_C }^*S_1,\dots ,\imath _{_C}^*S_q\}$ 
are basis for $H^0(C ,{\cal L}^l)$; 
\smallskip 

{\rm (iii)}  $S_{q+j}(x_j),S'_{q+j}(x_j)\neq 0$ for 
$j=1,\dots ,k$.
 
\nit 
The point $[(C ,\unl x,u)]\in \ovl M_{g,k}(X ,A)$ is $G$-stable 
if the same holds, but with strict inequalities in {\rm (i)}.
\end{proposition}

The shortcoming of this proposition is being too algebraic and 
therefore difficult to check in practice. For this reason we 
shall try to find necessary conditions on one hand and sufficient 
conditions on the other hand for it. 

Let us start with the necessary conditions. An easy consequence of 
the proposition above is the 

\begin{corollary}{\label{cor:plus-minus}} 
If the point $[(C,\unl x,u)]\in \ovl M_{g,k}(X,A)$ is 
$G$-semi-stable, then for all {\rm 1-PS} $\lambda :{\mbb C}^*
\rar G$ there are sections 
$$
S_j\in H^0(Y,{\cal L}^l)_{m_j}\kern1em j=1,\dots ,q
$$
satisfying the properties: 

\smallskip  
{\rm (i)} $\{ \imath_{_C }^*S_1,\dots ,\imath _{_C}^*S_q\} $ 
is a basis for $H^0(C ,{\cal L}^l)$;  

\smallskip 
{\rm (ii)} the set $\{ m_j\}_{j=1,\dots ,q}$ contains simultaneously 
positive and negative integers. 
\end{corollary}

\begin{proof} For a fixed 1-PS of $G$, there are two possibilities 
in the previous proposition: either all the $m_j$'s vanish for 
$j=1,\dots ,q$ and we are done, or it is not so. Assume that all 
$m_j\geq 0$ for $j=1,\dots ,q$. Because the sum 
$\sum_{j=1}^{j=q+k}m'_j\leq 0$, it follows that there must exist a 
$m'_{q+h}\leq 0$. We know that $S'_{q+h}(x_{q+h})\neq 0$ and 
therefore the restriction $\imath _{_C}^*S'_{q+h}\neq 0$. Because 
$\{ \imath_{_C }^*S_1,\dots ,\imath _{_C}^*S_q\} $ is 
a basis of $H^0(C ,{\cal L}^l)$, one can write  
$\imath _{_C}^*S'_{q+h}$ as a non-zero linear combination of 
these vectors. Now all we have to do is to replace a section 
from the set $\{S_1,\dots,S_q\}$ which appears in this linear 
combination with $S'_{q+h}$.
\end{proof}
 
The next proposition gives a geometrical restriction which must be 
satisfied by the $G$-semi-stable maps in $\ovl M_{g,k}(X,A)$.  

\begin{proposition}{\label{prop:Gss}}  
If $(C,\unl x,u)\in \ovl M_{g,k}(X,A)$ is a stable map which is
$G$-semi-stable, then for each {\rm 1-PS} $\lambda :{\mbb C}^*
\rar G$ there is an irreducible component $C_\delta$ of 
$C$ such that the image of the map 
$$
u_{|C_\delta }:C _\delta \lar X 
$$
is not contained in the $\lambda $-unstable locus of $X$. 
\end{proposition} 

\begin{proof} The line bundle ${\cal L}^l\rar Y $ is again very ample 
and its associated linear system gives an embedding 
$$
\begin{array}{ccl} 
{\cal L}^l &  & 
{\cal O}_{{\mbb P}^{R-1}}(1) \\[1ex]  
\dar & & \kern1em \dar \\ 
Y & 
\xrar{|{\cal L}^l|} & 
{\mbb P}^{R-1}={\mbb P}\left( {H^0(Y,{\cal L}^l)}^\vee \right) 
\end{array}
$$
The $G$-action on ${\cal L}^l\rar Y$ induces one on 
${\cal O}_{{\mbb P}^{R-1}}(1)\rar{\mbb P}^{R-1}$. For 
a stable map $(C,\unl x ,u)$ which is $G$-semi-stable and 
$\lambda :{\mbb C}^*\rar G$ a 1-PS of $G$, corollary 
\ref{cor:plus-minus} ensures the existence of sections 
$S_j\in H^0(Y,{\cal L}^l)_{m_j}$ whose restrictions to $C$ give a 
basis of $H^0(C,{\cal L}^l)$; in particular, they are linearly 
independent. Because of the direct sum decomposition 
$$
H^0(Y,{\cal L}^l)=\bigoplus _{m\in {\mbb Z}}H^0(Y,{\cal L}^l)_m,
$$
these sections can be completed with sections 
$$
S_{q+1}\in H^0(Y,{\cal L}^l)_{m_{q+1}}, \dots ,
S_R\in H^0(Y,{\cal L}^l)_{m_R}
$$
to a basis of $H^0(Y,{\cal L}^l)$. This basis defines coordinates 
on ${H^0(Y,{\cal L}^l)}^\vee \cong {\mbb C}^R$ in which the 
$\lambda $-action is diagonal.  
 
{\sl Claim } There exists an irreducible component $C_\delta$ 
of $C$ having the property that among $\{ S_1,\dots ,S_q\}$ 
there are two sections $S_j\in H^0(Y,{\cal L}^l)_{m_j}$ and 
$S_{j'}\in H^0(Y,{\cal L}^l)_{m_{j'}}$ such that $m_{j'}\leq 0$ and 
$m_j\geq 0$ and their restriction to $C_\delta$ is non-zero. 

We know already that there are two sections $S_\alpha$ and 
$S_\tau$ such that $m_\alpha \leq 0$ and $m_\tau \geq 0$, and 
their restriction to $C$ is nonzero. Let $C_\alpha$ and $C _\tau $ 
respectively two irreducible components of $C$ on which these two 
sections do not vanish. Because $C$ is connected, there is a chain 
of irreducible components $C_\alpha, C_\beta ,\dots , C_\tau $ 
connecting these two components. Since  $\{ S_1,\dots ,S_q\} $ is a 
basis of $H^0(C,{\cal L}^l)$ and ${\cal L}^l \rar C$ is very ample, 
it follows that there are sections $S_{\alpha\beta}, 
S_{\beta\gamma },\dots ,$ $S_{\sigma\tau}$ with the property that: 
$S_{\alpha \beta}$ does not vanish at a (certain) point in $C_\alpha 
\cap C_\beta$, $S_{\beta\gamma}$ does not vanish at a (certain) point 
in $C_\beta \cap C_\gamma ,\dots , S_{\sigma \tau }$ does not vanish 
at a (certain) point in $C _\sigma \cap C_\tau $. Notice that 
$S_{\alpha \beta}$ does not vanish on $C_\alpha $ and $C_\beta$, 
$S_{\beta\gamma}$ does not vanish on $C_\beta$ and $C_\gamma$ and so 
on. Let $m_{\alpha \beta}, m_{\beta \gamma },\dots $ denote the 
weights of the sections $S_{\alpha \beta }, S_{\beta \gamma },\dots $ 
respectively. If $m_{\alpha \beta} \geq 0$, then the component 
$C_\alpha $ satisfies the requirement of the claim. If it is not the 
case, we look at the chain $C _\beta , \dots , C_\tau $ whose length 
is one less than the length of $C_\alpha ,\dots ,C_\tau $. Because at 
the end $C_\tau $  the weight $m_\tau $ is positive, an induction 
argument on the length of the connecting chain shows that it must 
exist an irreducible component $C_\delta$ of the chain $C_\alpha ,
\dots , C_\tau $ having the property of the claim. 

Let us look now at the image of a point $p\in C_\delta $ inside 
${\mbb P}^{R-1}$: a representative $p'\in {\mbb C}^R $ of it will 
have non-zero coordinates with both positive and negative weights 
(for the $\lambda $ action), so $p$ is in the $\lambda $-semi-stable 
locus of $Y$. Since obviously $Y^\sst={\mbb P}^N\times X^\sst$, we 
deduce that $u(p)$ is in the $\lambda $-semi-stable locus of $X$. 
\end{proof}

In the case of a torus action, this proposition implies the  

\begin{corollary}{\label{cor:case-torus}} 
Suppose that a torus $T$ acts on $X$. If $(C,\unl x,u)\in \ovl 
M_{g,k}(X,A)$ is a $T$-semi-stable point and $C$ is irreducible, 
then the image of $u$ is not contained in the $T$-unstable locus 
of $X$. 
\end{corollary}

\begin{proof} By the Hilbert-Mumford criterion,
$$
X_G^\sst=\bigcap _{\lambda\ 1-PS\,{\rm of}\, T}X_\lambda ^\sst.
$$ 
Since $C$ is assumed irreducible, proposition~\ref{prop:Gss} implies 
that for any 1-PS $\lambda $ of $G$, the image of $u$ intersects the 
$\lambda $-unstable locus of $X$ in finitely many points; denote by 
$C^0(\lambda)$ the Zariski open subset of $C$ consisting of points 
which are mapped by $u$ into the $\lambda $-semi-stable locus of $X$. 
Because in a torus there are countably many one-parameter subgroups,
$$
C \neq \bigcup_{\lambda \ 1{\rm -PS}\,{\rm of}\, T}(C-C^0(\lambda )). 
$$
\end{proof}

In what follows we want to prove a weakened converse of 
proposition \ref{prop:Gss} which is useful when the
unstable locus $X^{\rm unstable}({\cal O}_X(1))$ has large
codimension in $X$. In this case it is reasonable to think 
that `many' curves in $X$ won't meet this locus at all. 

\begin{theorem}{\label{thm:ss-maps1}} 
A stable map $[(C,\unl x,u)]\in \ovl M_{g,k}(X,A)$ having the 
property that
$$
{\rm Image}\,(u:C\rar X)\subset X^\sst({\cal O}_X(1))
$$ 
is a $G$-semi-stable point of $\ovl M_{g,k}(X,A)$. 
\end{theorem} 

\begin{proof} The geometric invariant quotient $\wht X =X\invq G$ is 
a projective subvariety of $\wht {{\mbb P}^r}={\mbb P}^r\invq G$. 
This last geometric quotient can be described as 
$$
\wht {{\mbb P}^r}={\tt Proj}\left( 
\sum _n {H^0({\mbb P}^r,{\cal O}_{{\mbb P}^r}(n))}^G \right) .
$$
Let's denote by $\phi :{({\mbb P}^r)}^\sst\rar 
\wht{{\mbb P}^r}$ the quotient map. There is an invertible sheaf 
$M\rar \wht{{\mbb P}^r}$ such that $\phi ^*M
={\cal O}_{{\mbb P}^r}(m_0)_{|{({\mbb P}^r)}^\sst}$ for some $m_0>0$. 
It has the additional property that for large enough values of $n$,  
$$
H^0(\wht{{\mbb P}^r},M^n)\stackrel{\phi ^*}{\cong } 
{H^0({\mbb P}^r,{\cal O}_{{\mbb P}^r}(nm_0))}^G.
$$
The assumption that the image of the stable map is contained in the 
semi-stable locus of $X$ implies the existence of the commutative 
diagram 
$$
\begin{array}{rcl} 
(C,\unl x) & \stackrel{u}{\lar} & X^\sst\\  
_{\wht u} \sear\kern-2ex & &\kern-2ex\swar_\phi \\[1mm] 
 & \wht X. &
\end{array}
$$

\begin{remark}{\label{stab-quot-curves}} 
The map $\wht u:(C,\unl x)\rar \wht X$ is still stable. Indeed,
problems appear only if $\wht u$ contracts some 
${\mbb P}^1$-components, without enough special points on them, 
which are not contracted by $u$. If $C_\delta $ denotes such a 
${\mbb P}^1$-component of $C$,  
$$
{\rm deg}_{u_*C_\delta }{\cal O}_X(m_0)
={\rm deg}_{\wht u_*C_\delta }M=0,
$$
so that $u$ must be constant on $C_\delta$. This contradicts 
the stability of $u$.  
\end{remark} 

The group $G$ acts on $Y={\mbb P}^N\times {\mbb P}^r$ trivially on the 
first factor and consequently its invariant quotient is 
$\wht Y:={\mbb P}^N\times \wht{{\mbb P}^r}$. The quotient map 
$\psi :Y^\sst({\cal L})\longrightarrow \wht Y$ is just 
$\psi=({\rm id}_{{\mbb P}^N},\phi )$. Let us define the line bundle 
$$
{\cal M}:={\cal O}_{{\mbb P}^N}(m_0)\boxtimes M\lar \wht{Y}.
$$
It has the property that 
$$
\psi ^*{\cal M}={\cal O}_{{\mbb P}^N}(m_0)\boxtimes \phi ^*M
={\left( {\cal O}_{{\mbb P}^N}(m_0)\boxtimes 
{\cal O}_{{\mbb P}^r}(m_0)\right)}_{|Y^\sst}
={{\cal L}^{m_0}}_{|Y^\sst}
$$ 
and it can be easily checked that 
$$
H^0(\wht{Y},{\cal M}^n)\stackrel{\psi ^*}{\cong }
{H^0(Y,{\cal L}^{nm_0})}^G
$$
for large $n$. There is again a commutative diagram
$$
\xymatrix{
(C,\unl x)\ar[dr]_-{\hat u_Y}\ar[r]^-{u_Y}&
Y^\sst={\mbb P}^N\times {({\mbb P}^r)}^\sst
\ar[d]^{\psi=({\rm id}_{{\mbb P}^N},\phi)}\\
&\wht{Y}={\mbb P}^N\times \wht{{\mbb P}^r}.
}
$$
Because $u_{Y}$ is an embedding, $\wht u_{Y}$ is
also. The $1$-dimensional subvariety $\wht u_{Y*}C$ 
of $\wht Y$ has Hilbert polynomial: 
$$
\begin{array}{rl} 
\wht P(n) & = h^0(C,\wht u_{Y}^*{\cal M}^n)
-h^1(C,\wht u_{Y}^*{\cal M}^n)\\[1ex] 
 & = h^0(C,u_{Y}^*{\cal L}^{nm_0})
-h^1(C,u_{Y}^*{\cal L}^{nm_0})=P(nm_0),
\end{array}
$$
where $P$ is the Hilbert polynomial of $u_{Y*}C\subset Y$. It is 
independent of $(C,\unl x,u)\in \ovl M_{g,k}(X,A)$ satisfying the 
hypothesis of the proposition. 

Grothendieck proves in \cite{Gk} that there is an integer $k>0$ 
such that for all $n\geq k$, ${\cal M}^n$ is generated by its 
global sections and moreover, for any closed subscheme $\wht Z$ 
of $\wht Y$ whose Hilbert polynomial is $\wht P$ there is an 
epimorphism 
$$
H^0(\wht Y,{\cal M}^n)\longrightarrow  
H^0(\wht Y,{\cal O}_{\wht Z}\otimes {\cal M})\longrightarrow 0.
$$ 
Recall that for obtaining a projective embedding of 
$\ovl M_{g,k}({\mbb P}^r,d)$ we had to chose a high enough power 
${\cal L}^l\rar Y$. Since $\psi ^*{\cal M}
={{\cal L}^{m_0}}_{|Y^\sst}$, we can chose {\em from the very 
beginning} an integer $l$ large enough such that $\psi ^*{\cal M}^n
={{\cal L}^l}_{|Y^\sst}$ with $n\geq k$ ($l=nm_0$). 

The following three relations  
\begin{align}{\label{relations}}
\left\{ 
\begin{array}{l} 
H^0(C, u_{Y}^*{\cal L}^l)=H^0(C,\wht u_{Y}^*{\cal M}^n)
\\[1ex] 
H^0(\wht Y,{\cal M}^n)\lar H^0(C,\wht u_{Y}^*{\cal M}^n)\lar 0
\\[1ex] 
H^0(\wht Y,{\cal M}^n)\stackrel{\psi ^*}{\cong }
{H^0(Y,{\cal L}^l)}^G 
\end{array}
\right. 
\end{align}  
prove that there are  $G$-invariant sections $S_1, \dots ,S_q\in 
{H^0(Y,{\cal L}^l)}^G$ such that the restrictions 
$\{ \imath _{_C}^*S_1, \dots ,\imath _{_C}^*S_q\}$ form a basis of 
$H^0(C,u_{_Y}^*{\cal L}^l)$. The problem with the marked points is 
easy: by hypothesis $u(x_1),\dots ,u(x_k)\in X^\sst({\cal O}_X(1))$ 
and we may consider their images 
$\wht u_{_Y}(x_1),\dots,\wht u_{_Y}(x_k)\in\wht Y$. The 
number $n$ was chosen large enough to ensure that ${\cal M}^n $ is 
globally generated by its sections. Consequently, we find 
$\wht S_{q+1},\dots ,\wht S_{q+k}\in H^0(\wht Y,{\cal M}^n )$ such 
that $\wht S_{q+j}(x_j)\neq 0$ for $j=1,\dots, k$. Since 
$H^0(\wht Y,{\cal M}^n)\srel{\psi^*}{\cong}{H^0(Y,{\cal L}^l)}^G$, 
there are $G$-invariant sections 
$S_{q+1},\dots ,S_{q+k}\in {H^0(Y,{\cal L}^l)}^G$ such 
that $S_{q+j}(x_j)\neq 0$ for $j=1,\dots ,k$. 

The $q+k$ sections $S_1,\dots ,S_q,S_{q+1},\dots ,S_{q+k}$ now
obviously satisfy the conditions of the proposition 
\ref{prop:ss-maps}. 
\end{proof} 

\begin{corollary}{\label{cor:stable-maps}} 
If the stable map $[(C,\unl x,u)]$ has the property that 
${\rm Image}\,u\subset X_{(0)}^\s$, then 
$[(C,\unl x,u)]\in {\ovl M_{g,k}(X,A)}_{(0)}^\s$. 
\end{corollary} 

\begin{proof} 
There are two things to prove in this statement: the first one is 
that the stabilizer of $[(C,\unl x,u)]$ in $G$ is finite, and the 
second one is that this point is indeed stable.

When $k>0$, for any $x\in \unl x$ we have 
${\rm Stab}_G[(C,\unl x,u)]\subset{\rm Stab}_Gu(x)$, and therefore 
the stabilizer of the map is indeed finite. Let us prove that it is 
so in general. Consider a representative $u:(C,\unl x)\rar X$ of the 
point $[(C,\unl x,u)]$ and define $H:={\rm Stab}_G[(C,\unl x,u)].$  
Let's assume that $H$ is not finite. By definition, for any $h\in H$, 
there is an automorphism $\gamma _h\in {\rm Aut}(C,\unl x)$ having 
the property that $hu=u\gamma _h$. In particular, for all $h\in H$, 
Image$\,hu={\rm Image}\,u$. For $p\in C$ an arbitrary point, $u(p)$ 
has finite stabilizer in $G$ by assumption and therefore 
dim$\,H\!\cdot\!u(p)={\rm dim}\,H>0$. Since $H\!\cdot\! u(p)\subset 
{\rm Image}\,u$ which is one dimensional, we deduce that dim$\,H=1$. 
Let us look at the connected component of the identity $H^\circ $ 
of $H$: it is a connected $1$-dimensional group and therefore 
isomorphic either to the multiplicative group $G_m$ or to the 
additive group $G_a$. In both cases $\disp\lim_{t\rar\infty}
t\!\cdot\!u(p)\in {\rm Image}\,u$ will be fixed. This contradicts 
the assumption that Image$\,u\subset X_{(0)}^\s$. 

We will show that the point $[(C,\unl x,u)]$ is $G$-stable using 
proposition \ref{prop:ss-maps}. For a 1-PS $\lambda:{\mbb C}^*
\rar G$, there are sections $S_1\in H^0(Y,{\cal L}^l)_{m_1}$ and 
$S'_1\in H^0(Y,{\cal L}^l)_{m'_1}$ with $m'_1<0<m_1$, and whose 
restriction to $C$ is non-zero (this is because 
Image\,$u\subset X^\sst$). Using now \eqref{relations}, we 
complete $S_1$ and $S'_1$ with $G$-invariant sections 
$S_j,S'_j\in H^0(Y,{\cal L}^l)$ in order to fulfill the 
requirements of proposition \ref{prop:ss-maps}.
\end{proof} 


\section{The symplectic perspective of the problem}
{\label{sct:symp-pers}} 

In order to have a geometrically clearer picture of what is going 
on, we shall investigate the symplectic counterpart of the problem 
studied in the previous section. It is well-known that the invariant 
quotient in algebraic geometry has a very close analogue in 
symplectic geometry, namely the Marsden-Weinstein quotient. More 
precisely, assume that a complex, algebraic, linearly reductive group 
$G$ acts on a projective variety $X$ and the action is linearized  
in a very ample line bundle ${\cal O}_X(1)\rar X$. Let $K$ be the 
maximal compact subgroup of $G$ and denote $m:X\rar{\mfrak k}^*$ 
the moment map for the $K$-action, which takes values in the dual 
of the Lie algebra of $K$. Assume also that the $G$-action on $X$ is 
such that $X^\sst=X_{(0)}^\s$. Under these assumptions we have the 
\smallskip

{\sl Result } A point $x\in X$ is semi-stable if and only if 
$\ovl{Gx}\cap m^{-1}(0)\neq \emptyset $. Moreover, the inclusion 
$m^{-1}(0)/K\hookrightarrow X^\sst/G$ is a homeomorphism. 

Details and further references can be found in \cite{Ki}.
\smallskip 

It is clear that the map $\pi :\ovl M_{g,k}(X,A)\rar \ovl M_{g,k}$ 
which associates to a stable map $[(C,\unl x,u)]$ its stabilized 
curve $[(C^{st}, \unl x^{st})]$ is preserved by the $G$-action. 
In order to ensure that the space $\ovl M_{g,k}$ exists, we shall 
assume that $2g-2+k>0$. It is also clear that a point 
$[(C,\unl x ,u)]\in \ovl M_{g,k}(X,A)$ is $G$-semi-stable if and only 
if it is $G$-semi-stable when it is viewed as a point in the fibre 
$\ovl M_{g,k}(X,A)
\times_{\rm Spec\,{\mbb C}}k_{[(C^{st},\unl x^{st})]}$ of $\pi $, 
where $k_{[(C^{st},\unl x^{st})]}$ denotes the function field of 
the corresponding point. This remark justifies the following 
construction: for a quasi-stable curve $C$ of genus $g$, let us 
define 
$$
M_{C,k}(X,A):=\left\{   
u:(C,\unl x)\rar X \left| 
\begin{array}{l} 
(C,\unl x,u)\ {\rm is\ a\ stable\ map,}\\ 
|x|=k,\ u_*[C]=A  
\end{array}\right. \right\}. 
$$

\begin{lemma}{\label{lm:quasi-proj-global-cover}} 
$M_{C,k}(X,A)$ has a natural quasi-projective scheme structure. 
\end{lemma} 

\begin{proof} We may assume as usual that $X={\mbb P}^r$. Recall 
that a map $(C,\unl x,u)$ is stable if and only if $L_{(C,\unl x,u)}
=\omega_C(x_1+ \dots +x_k)\otimes u^*{\cal O}_{{\mbb P}^r}(3)
\rar C$ is ample. Also, there is an integer $f=f(g,k,r,d)>0$ 
with the property that $L_{(C,\unl x,u)}^f\rar C$ is very ample. In 
this way, any stable map $(C,\unl x,u)$ gave rise to an embedding 
$C\rar{\mbb P}\left( {H^0(C,L_{(C,\unl x,u)}^f)}^\vee\right)$ into 
a space {\em isomorphic} to ${\mbb P}^N$, where 
$N+1= {\rm dim}\,H^0(C,L_C^f)$. The ambiguity in the choice of this 
isomorphism is given by elements in $PGl(N+1)$. In order to define 
the space $M_{C,k}(X,A)$ we shall use a {\em fixed} but otherwise 
arbitrary stable map $(C,\unl x_0,u_0)$. Let us consider the 
{\em fixed embedding} 
\begin{align}{\label{fix-emb}}
C\srel{j_0}{\lar}
{\mbb P}\left({H^0\left(C,L_{(C,\unl x_0,u_0)}^f\right)}^\vee\right)
={\mbb P}^N
\end{align} 
defined by the very ample line bundle 
${\cal O}(1):=L_{(C,\unl x_0,u_0)}^f$ and let $e$ be its degree on 
$C$. For another stable map $(C,\unl x,u)$, the Hilbert polynomial 
of its graph $\Gamma_u\subset C\times {\mbb P}^r$ is 
$$
P(n)=\chi \left( {\cal O}_{\Gamma_u}\otimes 
{\cal O}_C(n)\boxtimes {\cal O}_{{\mbb P}^r}(n)\right) 
=n(e+d)+\chi (\Gamma _u)
=n(d+e)+(1-g).
$$
So each stable map $(C,\unl x,u)$ defines a point 
$(\Gamma_u,\unl x)\in {\cal H}ilb _{C\times {\mbb P}^r}^P
\times {(C\times {\mbb P}^r)}^k$. Using the embedding $j_0: 
C\rar {\mbb P}^N$, the graph $\Gamma _u$ 
can be viewed as a subvariety of ${\mbb P}^N\times {\mbb P}^r=Y$ 
and its Hilbert polynomial with respect to the very ample line 
bundle ${\cal L}={\cal O}_{{\mbb P}^N}(1)\boxtimes 
{\cal O}_{{\mbb P}^r}(1)\rar Y$ is $P$ also. 

Clearly, the same is true for any closed subscheme $Z$ of 
$C\times {\mbb P}^r$: the Hilbert polynomial inside 
$C\times {\mbb P}^r$ with respect to ${\cal O}_C(1)\boxtimes 
{\cal O}_{{\mbb P}^r}(1)$ is the same as the Hilbert polynomial of 
its image $j_{0*}Z$ inside $Y$ with respect to ${\cal L}$. So 
we obtain a closed immersion 
$$
{\cal H}_C:={\cal H}ilb_{C\times {\mbb P}^r}^P
\xrar{\phantom{MMM}}{\cal H}ilb_Y^P={\cal H}.
$$ 
The ample line on $\cal H$ is det$\,{\cal Q}\rar{\cal H}$, 
where $\cal Q$ is the universal quotient bundle on some Grassmann 
variety (see section~\ref{sct:reminders}). Remember that there is 
a subscheme $\cal S$ of ${\cal H}\times Y^k$ corresponding 
to the locus of $k$-marked, genus $g$ stable maps to ${\mbb P}^r$ 
which represent $d$ times the generator of $H_2({\mbb P}^r, 
{\mbb Z})$. Let consider the commuting diagram  
$$
\xymatrix{
\ovl M_{C,k}({\mbb P}^r,d)\kern-2.7em&
:={\cal H}_C\times{(C\times {\mbb P}^r)}^k
\times_{{\cal H}\times Y^k}{\cal S} 
\ar[d]\ar[r] & {\cal S}\ar[d] 
\\
&{\cal H}_C\times {(C\times {\mbb P}^r)}^k
\ar[r] & {\cal H}\times Y^k.
}
$$
Since $\ovl M_{C,k}({\mbb P}^r ,d)\rar {\cal H}_C
\times {(C\times {\mbb P}^r)}^k$ is an immersion, 
$\ovl M_{C,k}({\mbb P}^r, d)$ is a quasi-projective scheme. 
Its ample line bundle is determined by the restriction to 
${\cal S}$ of det$\,{\cal Q}_k\rar {\cal H}\times Y^k$. The 
space $M_{C,k}({\mbb P}^r ,d)$ is quasi-projective, being an 
open subset of $\ovl M_{C,k}({\mbb P}^r ,d)$. 
\end{proof}
 
\begin{remark} $\ovl M_{C,k}(X,A)$ is the Gromov compactification 
of $M_{C,k}(X,A)$. It is also projective if $(C,\unl x_0)$ is a 
{\em stable curve} in the sense of Deligne-Mumford (remember that 
we have used a fixed stable map $(C,\unl x_0,u_0)$ for the 
embedding $j_0:C\rar {\mbb P}^N$). In this case the intersection 
of the $PGl(N+1)$-orbit of a stable map in 
${\cal S}\rar {\cal  H}\times Y^k$ with the image of 
${\cal H}_C\times{(C\times {\mbb P}^r)}^k\rar {\cal H}\times Y^k$ 
consists of finitely many points. Consequently, the map 
$\ovl M_{C,k}({\mbb P}^r,d)\rar \ovl M_{g,k}({\mbb P}^r ,d)$ is 
generically finite on its image. Because 
$\ovl M_{g,k}({\mbb P}^r,d)$ is projective, 
the conclusion follows. 
\end{remark} 

The reason for introducing the space $M_{C,k}(X,A)$ is that now 
we can work with maps instead of equivalence classes of maps. 
For the symplectic point of view it is convenient to consider 
$M_{C,k}(X,A)$ with its reduced scheme structure, so that we look 
at it naively as being a quasi-projective variety. The Zariski 
tangent space of $M_{C,k}(X,A)$ at a point $(C,\unl x,u)$ is  
$$
T_{(C,\unl x,u)}M_{C,k}(X,A)\subset\left\{ 
(s,v_1,\dots ,v_k) \left| 
\begin{array}{l} 
s\in H^0(C,u^*TX )\\ 
v_j\in T_{x_j}C ,\ j=1,\dots ,k 
\end{array} \right. \right\} .  
$$
In the sequel we shall compute the K\"ahler form on 
$\ovl M_{C,k}(X,A)$ induced by its projective embedding. We can see 
that the ample line bundle det$\,{\cal Q}_k\rar{\cal H}\times Y^k$ 
introduced in section \ref{sct:reminders} and defined by 
\eqref{ample-line-bdl} is det$\,(p_*\ovl E^*{\cal L}^l)\otimes 
\ovl{\tt ev}^*({\cal L}^l)^{\boxtimes k}$, where  
$$
\begin{array}{ccl} 
\ovl E^*{\cal L}^l &  & 
{\cal L}^l ={\cal O}_{{\mbb P}^N}(l)\boxtimes 
{\cal O}_{{\mbb P}^r}(l)\\ 
\dar & & \kern3.5ex \dar \\ 
M_{C,k}(X,A)\times C& 
\xrar{\ovl E=(j_0,E)} & 
{\mbb P}^N\times {\mbb P}^r \\[1ex] 
p\dar \kern1ex & & \\[1ex]  
M_{C,k}(X,A) & & 
\end{array}
$$
and 
$$
M_{C,k}(X,A)
\xrar{\tt \ovl{ev} =(\ovl{ev}_1,\dots ,\ovl{ev}_k)}
{\left({\mbb P}^N\times {\mbb P}^r\right)}^k 
\kern1cm{\tt \ovl{ev}_j=({\it j}_0,ev _j)},
$$ 
$$
(C,\unl x,u)\lmt\left(
(j_0(x_1),u(x_1)),\dots,(j_0(x_k),u(x_k))
\right).
$$ 
The K\"ahler form on $M_{C,k}(X,A)$ induced by its projective
embedding is $-1/2\pi i \times [\text{curvature of 
det}((p_*\ovl E^*{\cal L}^l)\otimes \ovl{\tt ev}^*({\cal L}^l)
^{\boxtimes  k})]$. For computing this curvature, we need a 
Hermitian metric on ${\cal L}$ and a K\"ahler metric on the fibres 
of $M_{C,k}(X,A)\times C\stackrel{p}{\lar }M_{C,k}(X,A)$ i.e. on $C$. 
The fibres of $p$ will be all isometric, the K\"ahler form on them 
being 
\begin{align}{\label{metric-on-C}}
\gamma _{_C} :=\frac{1}{e}j_0^*\omega _{_{{\mbb P}^N}},
\; e:={\rm deg}_C{\cal O}_{{\mbb P}^N}(1).
\end{align}
This choice reflects the fact that for defining the space 
$M_{C,k}(X,A)$ we have required the maps $(C,\unl x,u)$ to have a 
{\em fixed} domain of definition. On ${\cal O}_{{\mbb P}^N}(1)$ 
and ${\cal O}_{{\mbb P}^r}(1)$ consider the Hermitian metrics whose 
curvatures are $-2\pi i\omega _{_{{\mbb P}^N}}$ and 
$-2\pi i \omega _{_{{\mbb P}^r}}$ respectively, with 
$\omega_{_{{\mbb P}^N}}$ and 
$\omega_{_{{\mbb P}^r}}$ the corresponding Fubini-Study forms. 
There is an induced Hermitian metric on 
${\cal L}={\cal O}_{{\mbb P}^N}(1)\boxtimes {\cal O}_{{\mbb P}^r}(1)$ 
and {\it a fortiori} on $\ovl E^*{\cal L}^l$. 

It is easy to see which is the expression of the curvature of 
$\ovl {\tt ev}^*{({\cal L}^l)}^{\boxtimes k}$ at a point 
$(C,\unl x,u)\in M_{C,k}(X,A)$: 
$$
\Omega _1:=-\frac{1}{2\pi i}
{( R^{{\tt ev}^*{({\cal L}^l)}^{\boxtimes k}})}_{(C,\unl x,u)}
=le\sum _{j=1}^k{(\gamma _{_C})}_{x_j}
+l\sum_{j=1}^k{\left({\tt ev}_j^*
\omega _{_{{\mbb P}^r}}\right)}_{(C,\unl x,u)}.
$$
For computing the curvature of det$(p_*\ovl E^*{\cal L}^l)$ the first 
thing to notice is that this line bundle is actually the determinant 
of the derived direct image of $\ovl E^*{\cal L}^l$, since by 
$p_*^1\ovl E ^*{\cal L}^l=0$ (see section~\ref{sct:reminders}). 
Consequently, we may apply the differential form of the Atiyah-Singer 
index theorem for families which is proved in a series of papers 
\cite{BF1,BF2,BGS} by J.-M.\,Bismut, D.\,Freed, H.\,Gillet, 
Ch.\,Soul\'e. According to \cite{BGS}, theorem 0.1 page 51, if $C$ is 
{\em a smooth curve}, the curvature 
$$
\Omega _2:=-\frac{1}{2\pi i}
R^{{\rm det}p_*\ovl E^*{\cal L}^l}
=\int _C Td\left( -\frac{1}{2\pi i}R^{T_{_C}}\right)  
\cdot 
exp\left( -\frac{1}{2\pi i}R^{\ovl E^*{\cal L}^l}\right).
$$
Here $R^{T_C}$ denotes the curvature 
of the relative tangent bundle of the projection $p$ (i.e. of $T_C$) 
corresponding to the K\"ahler metric $\gamma _{_C}$ on $C$ and 
$R^{\ovl E^*{\cal L}^l}$ is the curvature of the line bundle 
$\ovl E^*{\cal L}^l$ with respect to the Hermitian metric induced by 
that on ${\cal L}^l$. Let $\gamma :=(-1/2\pi i)R^{T_C}$; it is a real 
form of type $(1,1)$ on $C$ and therefore $\gamma =h\gamma _{_C}$ 
with $h:C\rar {\mbb R}$ a smooth function having the property 
that $\int _{_C}h\gamma _{_C}=2(1-g)$. On the other hand, 
$$
-\frac{1}{2\pi i}R^{\ovl E^*{\cal L}^l}
=l\ovl E^*\left( -\frac{1}{2\pi i}R^{\cal L}\right)
=l\ovl E^*(\omega _{_{{\mbb P}^N}}+\omega _{_{{\mbb P}^r}})
=l(e\gamma _{_C}+E^*\omega _{_{{\mbb P}^r}}).
$$
The form $\Omega _2$ is the degree two term in 
$$
\int _C\left( 1+\frac{1}{2}\gamma \right) 
\left( 1+l(e\gamma _{_C}+E^*\omega _{_{{\mbb P}^r}})
+\frac{l^2}{2}{(e\gamma _{_C}+E^*\omega _{_{{\mbb P}^r}})}^2\right).
$$
Making the computations we obtain:
$$ 
\begin{array}{rl}
\Omega _2 \hspace{-.7em} & 
\disp 
= \frac{l^2}{2}\int _C{(e\gamma _{_C}+E^*\omega _{_{{\mbb P}^r}})}^2
+\frac{l}{2}\int _C(e\gamma _{_C}+E^*\omega _{_{{\mbb P}^r}})\wedge 
\gamma \\[2ex] 
 & \disp 
=\frac{l^2}{2}\int _C{(E^*\omega _{_{{\mbb P}^r}})}^2 
+2eE^*\omega _{_{{\mbb P}^r}}\wedge \gamma _{_C} 
+\frac{l}{2}\int _CE^*\omega _{_{{\mbb P}^r}}\wedge \gamma \\[2ex] 
 & \disp 
=\frac{l^2}{2}\int _C{(E^*\omega _{_{{\mbb P}^r}})}^2 
+l^2e\int _C E^*\omega _{_{{\mbb P}^r}}\wedge \gamma _{_C} 
+\frac{l}{2}\int _C h\!\cdot\! E^*\omega _{_{{\mbb P}^r}}
\wedge\gamma _{_C}.
\end{array}
$$
This proves the 

\begin{proposition}{\label{prop:kahler-form-on-global-cover}} 
The curvature of the line bundle 
${\rm det}\,(p_*\ovl E^*{\cal L}^l)\otimes \ovl{\tt ev}^*
{({\cal L}^l)}^{\boxtimes k}$ divided by $-2\pi i$ is 
$\Omega =\Omega _1+\Omega _2$. It represents the K\"ahler form on 
$M_{C,k}(X,A)$ induced by the projective embedding 
described in lemma \ref{lm:quasi-proj-global-cover}. 
\end{proposition} 

Let us come back to the initial set-up: a complex, connected, 
li\-near algebraic group $G$ acts on a smooth, irreducible complex 
projective variety $X$ and the action is linearized in a very 
ample line bundle ${\cal O}_X(1)\rar X$. In this case, using the 
linear system associated to ${\cal O}_X(1)$, we may assume that 
$G$ acts on ${\mbb C}^{r+1}$ and, by an appropriate choice of 
coordinates, we are allowed to assume that the maximal compact 
subgroup $K$ of $G$ is included in $U(r+1)$. There is an induced 
action of $G$ (so, {\it a fortiori}, of $K$) on $M_{C,k}(X,A)$ 
defined by: $g\times (C,\unl x,u)\mapsto (C,\unl x,gu)$. Because the 
$K$-action preserves the Fubini-Study form $\omega _{_{{\mbb P}^r}}$ 
and the maps $\ovl E$ and $\ovl{\tt ev}$ are obviously both 
$K$-invariant, it follows that $\Omega $ is $K$-invariant also.  

We shall restrict ourselves to the case of ${\mbb C}^*$-actions; in a 
certain way, this is allowed by the Hilbert-Mumford criterion. The 
goal is to find a moment map for the induced $S^1$-action on 
$M_{C,k}(X,A)$. If such a moment map exists, it is uniquely defined,
up to a scalar constant by the K\"ahler form $\Omega $. 

The $S^1$-action on $X$ gives rise to a vector field $\xi$ on $X$ 
having the property that ${\cal L}_{\xi}\omega =0$ and 
${\cal L}_{\xi}J=0$ because $S^1$ acts on $X$ by isometries. The 
vector field $J{\xi}$ does not preserve $\omega $ in general but 
still preserves the complex structure of $X$. In fact the vector 
field $J\xi$ corresponds to the (holomorphic) action of 
${\mbb R}_+^*\hra {\mbb C}^*$. The  moment map $m:X\lar{\mbb R}$ 
corresponding to the action is $S^1$-invariant and  
$$
\rd m_x(J\xi_x)=\omega ({\xi}_x,J\xi_x)=\| {\xi}_x\|^2,
$$
so that $m$ is increasing along the flow lines of $J\xi$. 

Let us denote by $\cal V$ the vector field on $M_{C,k}(X,A)$ 
determined by the (holomorphic) $S^1$-action on $M_{C,k}(X,A)$. 
At a point $(C,\unl x,u)$,  
$$
{\cal V}_{(C,\unl x,u)}=(u^*{\xi},0,\dots ,0)
\in H^0(C,u^*TX)\times T_{x_1}C\times \dots T_{x_k}C.
$$ 
\begin{proposition}{\label{prop:mm-on-global-cover}} 
The function 
$$
\Psi:M_{C,k}(X,A)\lar {\mbb R},
$$ 
$$
\begin{array}{cl}
\Psi\kern-2mm& \disp 
:=l^2\int _C(m\circ E) E^*\omega _{_{{\mbb P}^r}} 
+l^2e\int _C(m\circ E)\gamma _{_C}\\[1ex]  
& \displaystyle \kern2mm  
+\frac{l}{2}\int_C h(m\circ E)\gamma _{_C} 
+l\sum _{j=1}^k m\circ{\tt ev}_j 
\end{array}
$$
is a moment map for this action. 
\end{proposition}

\begin{proof} Because $E,{\tt ev}_j,m,\omega _{_{{\mbb P}^r}}$ are 
all $S^1$-invariant, it follows that $\Psi$ is also. For proving 
that $\Psi$ is a moment map, we need to show that its 
differential is the same as the contraction of the K\"ahler 
form $\Omega $ on $M_{C,k}(X,A)$ with the vector field ${\cal V}$. 
It what follows, the symbol ``$\ort$'' will always denote the 
contraction of a differential form with a vector field. 

The contraction ${\cal V}\ort {(\gamma _{_C})}_{x_j}=0$ because the 
$T_{x_1}C\times \dots \times T_{x_k}C$-component of ${\cal V}$ is 
zero.\smallskip\nit 

$
\begin{array}{rl}
{\cal V}\ort ({\tt ev}_j^*\omega _{_{{\mbb P}^r}}) & 
={\tt ev}_j^*\left( ({\tt ev}_{j*}{\cal V})\ort 
\omega _{_{{\mbb P}^r}}\right) 
={\tt ev}_j^*({\xi}_{{\tt ev}_j(\cdot )}\ort 
\omega _{_{{\mbb P}^r}})\\[1ex] 
 & 
={\tt ev}_j^*({\rm d}m_{{\tt ev}_j(\cdot )})
={\rm d}(m\circ {\tt ev}_j). 
\end{array}
$\smallskip\nit     

$\disp 
{\cal V}\ort \int _C E^*\omega _{_{{\mbb P}^r}}\wedge \gamma _{_C} 
=\int _C{\cal V}\ort (E^*\omega _{_{{\mbb P}^r}}\wedge \gamma _{_C})
=\int _C\left( {\cal V}\ort E^*\omega _{_{{\mbb P}^r}}\right) \wedge 
\gamma _{_C}.$\smallskip
   
At a point $p\in C$,\smallskip\nit 

$
\begin{array}{rl} 
{\cal V}_p \ort {(E^*\omega _{_{{\mbb P}^r}})}_p & 
=E^*\left( E_*{\cal V}_p\ort 
{\omega _{_{{\mbb P}^r}}}_{,u(p)}\right) 
=E^*\left( {\xi}_{u(p)} \ort 
{\omega _{_{{\mbb P}^r}}}_{,u(p)}\right)\\[1.5ex] 
& 
=E^*\left( {\xi}_{u(p)} \ort 
{\omega_{_{{\mbb P}^r}}}_{,u(p)}\right)  
=E^*({\rm d}m_{u(p)})={\rm d}{(m\circ E)}_p.
\end{array}
$\smallskip
  
and consequently\smallskip\nit  

$
\begin{array}{cl}\disp  
{\cal V}\ort \int _CE^*\omega _{_{{\mbb P}^r}}\wedge \gamma _{_C}
\kern-2mm&\disp  
=\int _C{\rm d}(m\circ E)\wedge \gamma _{_C}
=\int _C{\rm d}\left( (m\circ E)\gamma _{_C}\right)\\[1ex] 
&\displaystyle  
={\rm d}\int_C(m\circ E)\gamma _{_C}.
\end{array}
$\smallskip\nit  

$ 
\begin{array}{rl}\disp 
{\cal V}\ort 
\int_C h\cdot E^*\omega _{_{{\mbb P}^r}}\wedge \gamma_{_C} & 
\disp 
=\int_C {\cal V}\ort 
\left( 
h\cdot E^*\omega _{_{{\mbb P}^r}}\wedge \gamma _{_C}
\right) 
\\[1ex] 
 &\kern-2cm\disp  
=\int_C h\cdot 
\left( {\cal V}\ort E^*\omega _{_{{\mbb P}^r}}\right) \wedge 
\gamma _{_C} 
=\int_C h\cdot {\rm d}(m\circ E) \wedge \gamma _{_C}
\\[1.5ex] 
&\kern-2cm\disp 
\stackrel{(\star)}{=}
\int _C{\rm d}\left( h(m\circ E)\gamma _{_C}\right) 
={\rm d}\int _C h(m\circ E)\gamma _{_C}. 
\end{array}
$\smallskip

For writing equality $(\star )$, we have used that 
${\rm d}h\wedge\gamma_{_C}={\rm d}_{_C}h\wedge\gamma_{_C}=0.$ 
\smallskip\nit 

$
\begin{array}{rl}\disp 
{\cal V}\ort \int _C{\left( E^*\omega _{_{{\mbb P}^r}}\right)}^2 & 
\disp 
=2\int _C\left( {\cal V}\ort E^*\omega _{_{{\mbb P}^r}}\right) \wedge 
E^*\omega _{_{{\mbb P}^r}} 
=2\int _C{\rm d}(m\circ E)E^*\omega _{_{{\mbb P}^r}} \\[2ex] 
 & \disp
=2\int _C{\rm d}\left( (m\circ E)E^*\omega _{_{{\mbb P}^r}}\right) 
={\rm d}\left( 2\int _C(m\circ E)E^*\omega _{_{{\mbb P}^r}} \right) 
\end{array}
$\smallskip

All together, these equalities show that $\Psi$ is indeed a moment 
map.
\end{proof}

\begin{remark}{\label{rk:limit-mm}}
It is interesting to look at the limit moment map for large $e$ and 
$l$. We should recall that $e$ is defined in \eqref{fix-emb} as the 
degree of a certain very ample line bundle on $C$ used to get a fixed 
embedding of $C$ into a projective space, while $l$ is an integer 
large enough for obtaining a projective embedding of the Hilbert 
scheme $\cal H$ (actually $l$ does depend on $e$). Because $X$ is 
compact, the moment map $m$ is bounded and therefore the last two 
terms in $\frac{1}{el^2}\Psi$ are of order 
$O\bigl(\frac{1}{el}\bigr)$. We will show that the first term is of 
order $O\bigl(\frac{1}{e}\bigr)$. Indeed, since maps 
$u\in M_{C,k}(X,A)$ are holomorphic, the pull-back 
$u^*\omega_{_{{\mbb P}^r}}\geq 0$ as a form on $C$. This implies that 
$$
\left|\frac{1}{e}\int_C (m\circ u)u^*\omega_{_{{\mbb P}^r}}\right|
\leq \frac{\max _{_X}|m|}{e}\int_Cu^*\omega_{_{{\mbb P}^r}}
=\frac{d\max _{_X}|m|}{e}=O\biggl(\frac{1}{e}\biggr),
$$
where $d$ denotes, as usual, the degree of the composite map 
$C\rar X\rar{\mbb P}^r$, which is a constant. 

The conclusion of this discussion is that for large $e$ and $l$ 
\begin{align}{\label{eq:limit-mm}}
\frac{1}{el^2}\Psi\sim\int _C(m\circ E)\gamma_{_C},
\end{align}
so the zero set of $\Psi$ will be close to the zeros of this second 
function. Notice that the right hand side of \eqref{eq:limit-mm} is 
the moment map corresponding to the K\"ahler form 
$$
\Omega_{\infty}:=\int_C E^*\omega_{_{{\mbb P}^r}}
\wedge\gamma_{_C}
$$ 
on $M_{C,k}(X,A)$.

At this point some care is required because $\gamma_{_C}$, as it 
is defined by \eqref{metric-on-C}, does depend on $e$. However, it 
is well-known (see \cite{Tian}) that the sequence of such metrics 
converges to a metric on $C$, which was still denoted $\gamma_{_C}$. 
\end{remark}

\begin{lemma}{\label{lm:finite-stab-on-global-cover}} 
$S^1$ acts with finite stabilizers on $\Psi^{-1}(0)$. 
\end{lemma}

\begin{proof} We claim that for $(C,\unl x,u)\in\Psi^{-1}(0)$ 
there is a point $p\in C$ with the property that $m(u(p))=0$ i.e. 
$u(p)\in m^{-1}(0)$. Suppose that it is not the case, so either 
Image$\,u\subset \{ m<0\} $ or Image$\,u\subset \{ m>0\} $. Let 
assume that we are in the first case. At a point 
$(C,\unl x,u)\in M_{C,k}(X,A)$, 
$$
\Psi(C,\unl x,u)=
l^2\kern-3pt\int _C\kern-1pt (m\circ u)u^*\omega _{_{{\mbb P}^r}}
+l\kern-3pt\int_C\kern-3pt
\left( le+\frac{1}{2}h\right) 
(m\circ u) \gamma _{_C}
+l\sum _{j=1}^km(u(x_j)).
$$
Recall that the smooth real-valued function $h$ defined on 
$C$ is the ``quotient'' $R^{T_{_C}}/\gamma_{_C}$, where 
$R^{T_{_C}}$ denotes the curvature of the tangent bundle of $C$ 
with respect to the K\"ahler form $\gamma _{_C}$. This last form 
was defined in terms of a {\em fixed} projective embedding of $C$; 
in particular, it does not depend on $l$. For obtaining the 
projective embedding of the Hilbert scheme we had to take large 
positive integral values for $l$ and therefore we may assume that 
$l$ is large enough for $le+\frac{1}{2}h$ to be a strictly positive 
function on $C$. Notice also that since $u$ is holomorphic and 
$\omega _{_{{\mbb P}^r}}$ is a positive $(1,1)$-form, the 
$(1,1)$-form $u^*\omega _{_{{\mbb P}^r}}$ on $C$ is still positive. 
It becomes now clear that if Image$\,u\subset \{ m<0\} $, 
$\Psi(C,\unl x,u)$ will be negative also. This contradicts the choice 
of $(C,\unl x,u)$ in the zero locus of $\Psi$.  

The lemma follows now because, by assumption, $S^1$ acts with finite 
stabilizers on $m^{-1}(0)$. 
\end{proof} 

Using these symplectic techniques, we recover easily the 

\begin{theorem}[\ref{thm:ss-maps1}]{\label{thm:geom-picture}} 
A stable map $(C,\unl x,u)$, with $C$ smooth, having the 
property that 
${\rm Image}\,u\subset X_{(0)}^\s({\cal O}_X(1))$, 
defines a $G$-stable point in $\ovl M_{g,k}(X,A)$.
\end{theorem} 

\begin{proof} According to the Hilbert-Mumford criterion, 
is sufficient to prove the statement for every 1-PS 
$\lambda :{\mbb C}^*\rar G$. For a fixed 1-PS 
$\lambda $ of $G$, the point $(C,\unl x,u)$ is $\lambda$-stable, 
if its ${\mbb C}^*$-orbit meets the zero-level set of the moment 
map $\Psi$ on $M_{C,k}(X,A)$. Assume, for instance, that 
$\Psi(C,\unl x,u)<0$. By hypothesis, Image$\,u\subset X_{(0)}^\s
\subset X^{\s(\lambda )}$, so that under the ${\mbb R}_+^*$-action 
all the points $u(p),\,p\in C,$ meet the $m^{-1}(0)$-level; 
consequently $m(r\cdot u(p))>0$ for all $p\in C$ and $r\gg 0$. 
For such a large $r$, the translated map $ru$ will have the 
property that $\Psi(C,\unl x,ru)>0$. A continuity argument proves 
that there is (a unique) $r_0$ such that 
$(C,\unl x,r_0u)\in\Psi^{-1}(0)$. Now, by 
lemma~\ref{lm:finite-stab-on-global-cover}, 
$(C,\unl x,u)$ has also finite stabilizer.
\end{proof}

\section{First application: comparison of invariants}
{\label{sct:conseq1}}

In this section we shall use the results obtained so far, comparing 
the genus zero Gromov-Witten invariants of a projective manifold with 
those of its invariant quotient for a group action. The main result is 
the following

\begin{theorem}{\label{thm:comparison}}
Consider a complex, connected, linearly reductive group $G$ 
acting on the irreducible projective variety $X$, and also a 
linearization of the action in a very ample line bundle 
${\cal O}_X(1)\rar X$. Denote by $\zeta\in H^{2\dim G}(X;{\mbb Q})$ 
the class of a rational transverse slice to $X\dashto \hat X$. 
Let $A\in H_2(X;{\mbb Z})$ be a class which can be represented by a 
morphism ${\mbb P}^1\rar X^\sst$ and denote $\hat A$ the class of 
its image in $\hat X$.

Suppose that the following conditions are satisfied:

{\rm (a1)\,} $G$ acts freely on the $G$-semistable locus of 
$X$, so that the quotient map $X^\sst\rar X\invq G=:\hat X$ 
is a principal $G$-bundle;

{\rm (a2)\,} $\ovl M_{0,k}(X,A)$ is generically smooth and has 
the expected dimension;

{\rm (a3)\,} every irreducible component of $\ovl M_{0,k}(X,A)$ 
contains a point represented by a morphism ${\mbb P}^1\rar X^\sst$;

{\rm (a4)\,} 
$M_{0,k}(\hat X,\hat A)\subset\ovl M_{0,k}(\hat X,\hat A)$ 
is a dense open subset.

\nit Then for any $\hat\alpha\in H^*(\hat X^k;{\mbb Q})$ the 
following equality between the genus zero Gromov-Witten invariants 
holds:
$$
GW_{\hat X,\hat A}^{0,k}(\hat\alpha)
=GW_{X,A}^{0,k}(\alpha\cup(pr^{X^k}_X)^*\zeta),
$$
where $\alpha\in H^*(X^k;{\mbb Q})$ is obtained from 
$\hat\alpha$ using the correspondence induced by the rational 
map $X^k\dashto\hat X^k$, and $pr^{X^k}_X:X^k\rar X$ denotes 
the projection onto the first component.
\end{theorem}

\nit Before proceeding to the proof of the theorem, I would like 
to discuss\smallskip

{\sl When are the hypothesis in theorem 
\ref{thm:comparison} satisfied?\,}  
The condition (a1) on $G$ to act freely on the semi-stable locus 
of $X$ is necessary in order to ensure the equality of the expected 
dimensions of the spaces of stable maps involved inhere. It is 
unlikely to have any relations between the invariants 
if there are semi-stable points in $X$ with positive dimensional 
stabilizers. I have imposed the condition (a2) in order to avoid the 
use of the virtual class, which could be a rather difficult task in 
the present context. Condition (a3) excludes the existence of 
irreducible components of $\ovl M_{0,k}(X,A)$ such that the images 
of {\em all} the corresponding morphisms cut the unstable locus 
of $X$. Condition (a4) says that there should be no irreducible 
component of $\ovl M_{0,k}(\hat X,\hat A)$ such that {\em all} its 
points represent stable maps whose domain of definition are trees 
of ${\mbb P}^1$'s. 

I would like now to enumerate some cases where the theorem above 
applies.

\begin{lemma}{\label{lm:conv}} 
The hypothesis {\rm (a2)-(a4)} in \ref{thm:comparison} 
are satisfied in the following cases: 

{\rm (i)\,} $\ovl M_{0,k}(X,A)$ has expected dimension and 
both of $\ovl M_{0,k}(X,A)$ and $\ovl M_{0,k}(\hat X,\hat A)$ 
are irreducible or 

{\rm  (ii)\,} when $X$ and $\hat X$ are both convex varieties; 
this is the case when either:
 
\quad {\rm (iia)} $X$ is convex and $T_{X^\sst}$ is generated by 
$G$-invariant sections or 

\quad {\rm (iib)} $X$ is convex and $G$ is simply connected.
\end{lemma}

\begin{proof}
(i) the conditions (a2)-(a4) follow immediately from the 
irreducibility of $\ovl M_{0,k}(X,A)$ and from the initial 
assumption that the class $A$ is can be represented by a 
morphism ${\mbb P}^1\rar X^\sst$. 

(ii) If both $X$ and $\hat X$ are convex, 
$M_{0,k}(X,A)\subset\ovl M_{0,k}(X,A)$ and 
$M_{0,k}(\hat X,\hat A)\subset\ovl M_{0,k}(\hat X,\hat A)$ 
are open and dense (see \cite{FP}, theorem 2, page 56); 
convexity implies also that we are working in the expected 
dimension. The only thing to check is condition (a3): this 
follows from the fact that the evaluation map at the 
$(k+1)^{\text{th}}$ marked point on $M_{0,k+1}(X,A)$ is 
submersive, and therefore any map ${\mbb P}^1\rar X$ can be 
`pulled away' from the unstable locus of $X$.

For both (iia) and (iib) we must prove that $\hat X$ is 
still convex. In the first case, we use the exact sequence 
$$
0\lar {\cal O}({\cal L}ie\,G)
\lar T_{X^\sst}^{\text{inv}}
\lar T_{\hat X}\lar 0
$$
on $\hat X$ associated to the principal bundle 
$X^\sst\rar\hat X$. In the second case we notice that for a 
morphism $v:{\mbb P}^1\rar\hat X$, $v^*X^\sst\rar{\mbb P}^1$ 
is a principal $G$-bundle which is topologically trivial since 
$G$ is simply connected; a result due to Grothendiek says 
that in this case $v^*X^\sst\rar{\mbb P}^1$ is in fact 
holomorphically trivial. This implies that the morphism 
$v$ can be lifted to a morphism $u:{\mbb P}^1\rar X^\sst$. 
The conclusion follows now from the convexity of $X$ and 
the exact sequence 
\begin{align}{\label{eqn:exact}}
0\lar T_{\mfrak g}
\lar T_{X^\sst}
\lar \phi^* T_{\hat X}\lar 0
\end{align}
on $X$, where $T_{\mfrak g}$ denoted the trivial sub-bundle 
generated by the infinitesimal action of $G$.
\end{proof}

Let us return now to the proof of theorem \ref{thm:comparison}, 
and consider the following commutative diagram: 
\begin{align}{\label{key}}
\begin{array}{ccccc}
\ovl M_{0,k}(X,A) & 
\xrar{\phantom{M}{\tt ev}_k^X\phantom{M}} & 
V & 
\xrar{\phantom{M}{\jmath _{V}}\phantom{M}} & 
X^k
\\ 
\kern-1ex\text{\scriptsize\rmR}'\vdashto &  & 
\vdashto &  &\kern0.5ex\vdashto{\rmr '}
\\ 
\wht{\ovl M_{0,k}(X,A)} & 
\srel{\wht{\kern2pt{{\tt ev}_k^X}\kern2pt}}{\longdashar} & 
\wht V & 
\xrar{\phantom{M}\jmath _{{\wht V}}\phantom{M}} & 
\wht{\kern2pt X^k}
\\ 
\kern-1.5ex\text{\scriptsize\rmR}''\vdashto & & & &
\kern1.5ex \vdashto{\rmr ''}
\\ 
\ovl M_{0,k}(\wht X,\wht A) & 
\xrar{\phantom{M}{\tt ev}_k^{\wht X}\phantom{M}}& 
W & 
\xrar{\phantom{M}{\jmath _{W}}\phantom{M}} & 
\kern1ex{\hat X}^k.
\end{array}
\end{align}
The notations are as follows: $V$ and $W$ are respectively the 
images of the morphisms ${\tt ev}_k^X$ and ${\tt ev}_k^{\wht X}$, 
both with the reduced scheme structure. The group $G$ acts on 
$X^k$ in a diagonal fashion and the eva\-luation morphism 
${\tt ev}_k^X$ is $G$-equivariant. The invariant quotients of 
$\ovl M_{g,k}(X,A)$, $V$ and $X^k$ are denoted respectively 
$\wht{\ovl M_{g,k}(X,A)}$, $\wht V$ and $\wht{\kern2pt X^k}$. 
Notice that $\wht V\neq \emptyset$ as soon as there are stable 
maps whose image is contained in the $G$-semi-stable locus of $X$. 
The universality property of quotients implies the existence of 
the rational map $\wht{\kern2pt{\tt ev}_k^X\kern2pt}$. The quotient 
map $X^k\dashto\hat X^k$ naturally factorizes through a rational map 
${\rmr ''}:\wht{\kern2pt X^k}\dashto\hat X^k$ whose general fibre is 
isomorphic to $G^k/G$. If $G$ was a torus, then ${\rmr ''}$ would 
have been the quotient map for the induced $G^k/G$-action on 
$\wht{\kern2pt X^k}$.

\begin{proposition}{\label{prop:birational}}
Under the assumptions {\rm (a1)-(a4)} in \ref{thm:comparison}, 
the map 
$$
\rmR'':\ovl M_{0,k}(X,A)\invq G
\dashto\ovl M_{0,k}(\hat X,\hat A)
$$ 
is birational.
\end{proposition}

\begin{proof}
I start noticing that $\rmR''$ is generically injective on its image. 
Indeed, by assumption (a3), we may restrict our attention to maps 
whose image is contained in $X^\sst$: if $u_1,u_2:{\mbb P}^1\rar 
X^\sst$ are morphisms such that $\phi\circ u_1=\phi\circ u_2$, it 
follows from (a1) that there is a morphism $g:{\mbb P}^1\rar G$ 
such that $u_2(\zeta)=g(\zeta)u_1(\zeta)$ for all 
$\zeta\in{\mbb P}^1$. Since $G$ is affine, the morphism $g$ must 
be constant and therefore $u_1$ and $u_2$ represent the same 
point in $\ovl M_{0,k}(X,A)\invq G$.

Write
$$
\ovl M_{0,k}(X,A)=\bigcup_\nu \ovl M_{0,k}(X,A)_\nu
$$
as the union of its irreducible components.  Assumption (a3) says 
that each component $\ovl M_{0,k}(X,A)_\nu$ contains a non-empty open 
subset $\ovl M_{0,k}(X,A)_\nu^o$ having the property that its points 
represent stable maps defined on ${\mbb P}^1$, with image completely 
contained in the $G$-stable locus of $X$. 

We know already that the maps 
$\rmR''_\nu:M_{0,k}(X,A)_\nu^o\rar \ovl M_{0,k}(\hat X,\hat A)$ 
are birational on their image; let us denote 
$\ovl M_{0,k}(\hat X,\hat A)_\nu$ the closures of these images. 
They are {\em distinct} irreducible components of 
$\ovl M_{0,k}(\hat X,\hat A)$, 
each of them having expected dimension. In fact, for 
$[({\mbb P}^1,\unl x,u)]\in\ovl M_{0,k}(X,A)_\nu^o$ the composite 
$[({\mbb P}^1,\unl x,\hat u)]\in\ovl M_{0,k}(\hat X,\hat A)_\nu$ 
is a smooth point of $\ovl M_{0,k}(\hat X,\hat A)$; this can be 
seen pulling back by $u$ the exact sequence \eqref{eqn:exact}.

For proving that $\rmR''$ is dominant, we have to show that if 
$\ovl M_{0,k}(\hat X,\hat A)_\mu$ is an irreducible component  
of $\ovl M_{0,k}(\hat X,\hat A)$, then $\mu$ is one of the 
$\nu$'s comming from $X$. This is the place where we are 
using the assumption (a4) which says that we may restrict 
ourselves to morphisms $\hat u:{\mbb P}^1\rar \hat X$. The 
{\em topological} type of $\hat u^*X^\sst\rar {\mbb P}^1$ is 
determined by the class $\hat u_*[{\mbb P}^1]=\hat A=\phi_*A$, 
and we assumed that $A$ can be represented by a morphism 
${\mbb P}^1\rar X^\sst$. Consequently, for any morphism 
$\hat u:{\mbb P}^1\rar \hat X$ representing the class $\hat A$, 
the principal bundle $\hat u^*X^\sst$ is topologically trivial 
and therefore holomorphically trivial by the same result of 
Grothendieck.
\end{proof}

\begin{definition}{\label{def:corresp}}
Let $f:M\dashto N$ be a rational map. We define 
$f^*:H^*(N)\rar H^*(M)$ by 
$f^*\alpha:={\rm PD}_M\bigl( p_*(\Gamma_f\cap q^*\alpha)\bigr)$, 
where $\Gamma_f\subset M\times N$ denotes the closure of the graph 
of $f$ with the projections $p$ and $q$ on $M$ and $N$ respectively 
and {\rm PD} stands for Poincar\'e-duality.
\end{definition}

This is just another way to express the correspondence induced by 
$f$. We should keep in mind that $f^*$ is not a ring homomorphism 
in general. 

For the proof of theorem \ref{thm:comparison}, we remark that if 
$V\subset M$ and $W\subset N$ are complete, irreducible subvarieties 
such that $V\cap{\cal D}om(f)\neq\emptyset$ and the restriction 
$f_V:V\dashto W$ is dominant, then for any $\alpha\in H^*(N)$ 
$$
\langle f^*\alpha,[V]\rangle
=\deg(f_V)\langle \alpha,[W]\rangle,
$$
where we set $\deg(f_V)=0$ when $\dim V>\dim W$. 
This claim follows from the fact that 
$(V\times N)\cdot \Gamma_f=\Gamma_{f_V}$. 

\begin{proof}(of theorem \ref{thm:comparison})\; 
We have seen in proposition \ref{prop:birational} that 
$$
\ovl M_{0,k}(X,A)
=\bigcup_\nu \ovl M_{0,k}(X,A)_\nu,\quad 
\ovl M_{0,k}(X,A)\invq G
=\bigcup_\nu \ovl M_{0,k}(X,A)_\nu\invq G
$$
and 
$$
\ovl M_{0,k}(\hat X,\hat A)
=\bigcup_\nu \ovl M_{0,k}(\hat X,\hat A)_\nu.
$$
Moreover, $\rmR ''_\nu:\ovl M_{0,k}(X,A)_\nu\invq G
\dashto \ovl M_{0,k}(\hat X,\hat A)_\nu$ are birational for all 
$\nu$. Let $V_\nu\subset X^k$ and $W_\nu\subset \hat X^k$ be 
respectively the images of the $k$-point evaluation maps on 
$\ovl M_{0,k}(X,A)_\nu$ and $\ovl M_{0,k}(\hat X,\hat A)_\nu$. 
Then $V_\nu\invq G$ is the closure of the image of the evaluation 
map on $\ovl M_{0,k}(X,A)_\nu\invq G$ and there is a natural map 
$\rmr ''_\nu:V_\nu\invq G\dashto W_\nu$ compatible with the other 
arrows in \eqref{key} which is dominant. 

The class $\zeta$ which appears in the statement of the theorem 
is just the class of a ``rational section'' of the quotient 
$\phi$, that is $\zeta:=\frac{1}{d}Z$ for a general complete 
intersection $Z\hra X$ which transversally intersects, in 
$d$ points, the closures of the general $G$-orbits in $X$. With 
this choice for $Z$, the rational map 
$V_\nu\cap(Z\times X^{k-1})\dashto V_\nu\invq G$ 
is generically finite of degree $d$. 

Consider $\hat\alpha\in H^*(\hat X^k)$ and let 
$\alpha:=\phi^*\hat\alpha$ for $\phi:X^k\dashto\hat X^k$. The 
discussion preceding this proof applied to the composite 
$V_\nu\cap(Z\times X^{k-1})\srel{d:1}{\dashto}
V_\nu\invq G\srel{\rmr ''_\nu}{\dashto}W_\nu$ shows that 
$$
\langle\alpha,[V_\nu\cap(Z\times X^{k-1})]\rangle
=d\cdot\deg(\rmr '')\langle\hat\alpha,[W_\nu]\rangle.
$$
We can write therefore,
$$
\begin{array}{cl}
\disp\int_{\ovl M_{0,k}(\hat X,\hat A)_\nu}\kern-3ex
{({\tt ev}_k^{\hat X})}^*ß\hat\alpha 
\kern-1ex
&\disp
=\deg({\tt ev}_k^{\hat X})\langle\hat\alpha,[W_\nu]\rangle
=\frac{\deg({\tt ev}_k^{\hat X})}{d\cdot\deg(\rmr ''_\nu)}
\langle\alpha,[V_\nu\cap(Z\times X^{k-1})]\rangle
\\[1ex]
&\kern-2.5em
=\deg(\wht{\kern2pt{{\tt ev}_{k,\nu}^X}\kern2pt})
\langle\alpha,[V_\nu]\cap{({\rm pr}_X^{X^k})}^*\zeta\rangle
=\deg({\tt ev}_{k,\nu}^X)
\langle\alpha,[V_\nu]\cap{({\rm pr}_X^{X^k})}^*\zeta\rangle
\\[1ex]
&\kern-2.5em
\disp =\deg({\tt ev}_{k,\nu}^X)
\langle\alpha\cup{({\rm pr}_X^{X^k})}^*\zeta,[V_\nu]\rangle
=\int_{\ovl M_{0,k}(X,A)_\nu}\kern-3ex
{({\tt ev}_k^X)}^*(\alpha\cup{({\rm pr}_X^{X^k})}^*\zeta).
\end{array}
$$
This finishes the case when $\deg(\rmr ''_\nu)\neq 0$. When 
$\deg(\rmr ''_\nu)=0$ (that is $\rmr ''$ is not generically 
finite), both sides are zero. Summing these equalities after 
$\nu$ we get the conclusion. 
\end{proof} 

\subsection{\sl Some examples}
{\label{s-sct:expl}}

The requirements in theorem \ref{thm:comparison} make its 
applications rather resticted. I shall present below some cases 
when all the hypothesis are fulfilled. 

\subsubsection{}
{\label{ss-sct:proj-space}} 

For $m,n\geq 1$ two integers, consider the linearized 
${\mbb C}^*$-action on ${\mbb P}^{m+n+1}$ given by 
$$
{\mbb C}^*\times {\mbb C}^{m+n+2}\lar {\mbb C}^{m+n+2},
$$ 
$$
t\times (z,w)=(tz,t^{-1}w)\quad 
\forall z\in {\mbb C}^{m+1},\forall w\in {\mbb C}^{n+1}.
$$
The unstable locus of ${\mbb P}^{m+n+1}$ is the union of the  
planes $L':=\{z=0\}$ and $L'':=\{w=0\}$ of codimension $m+1$ and 
$n+1$ respectively. The stabilizer of any semi-stable point is 
$\{+1,-1\}$, so the semi-stable locus coincides with the properly 
stable one. The (geometric) quotient of ${\mbb P}^{m+n+1}$ is 
${\mbb P}^m\times{\mbb P}^n$ and the quotient map is 
$$
\phi :{({\mbb P}^{m+n+1})}^\s\lar {\mbb P}^m\times{\mbb P}^n,
$$ 
$$
\phi ([z,w])=[z]\times [w].
$$
Any stable map $u:{\mbb P}^1\rar {({\mbb P}^{m+n+1})}^\s$ 
of degree $d$ induces a stable map of bidegree $(d,d)$ 
into the quotient. 

Consider the simplest case when $m=n=1$. Say that we look at the 
rational curves of degree $d$ in ${\mbb P}^3$ and {\it a fortiori} 
of bidegree $(d,d)$ in its quotient which is ${\mbb P}^1\times 
{\mbb P}^1$. The right number of marked points is $k=4d-1$. The 
class of a point in ${\mbb P}^1\times {\mbb P}^1$ is mapped by 
the correspondence 
$A^*({\mbb P}^1\times {\mbb P}^1)\vdash A^*({\mbb P}^3)$ into the 
class of a line, because a hyperplane intersects the general 
${\mbb C}^*$-orbit in ${\mbb P}^3$ in only one point. The class 
$\zeta $ is just the class of the hyperplane in ${\mbb P}^3$. 
According to the result, 
$$
GW_{{\mbb P}^3,d}^{4d-1}({\rm point}\times\underbrace{{\rm line}
\times\dots\times {\rm line}}_{4d-2\ {\rm times}})
=
GW_{{\mbb P}^1\times{\mbb P}^1,(d,d)}^{4d-1}
(\underbrace{{\rm point}\times\dots\times
{\rm point}}_{4d-1\ {\rm times}}). 
$$ 
In more down-to-earth terms, the number of rational curves of 
bidegree $(d,d)$ in ${\mbb P}^1\times {\mbb P}^1$ passing through 
$4d-1$ general points equals the number of rational curves of 
degree $d$ in ${\mbb P}^3$ passing through a point and another 
$4d-2$ general lines. 

\subsubsection{}
{\label{ss-sct:gras}}  

Given two integers $m>n>0$, consider the linearized 
$Sl_n({\mbb C})$-action on ${\mbb P}^{mn-1}:=
{\mbb P}\left({\Hom}({\mbb C}^m,{\mbb C}^n)\right)$ given by 
$$
Sl_n({\mbb C})\times \Hom({\mbb C}^m,{\mbb C}^n)
\lar \Hom({\mbb C}^m,{\mbb C}^n),\; 
(g,A)\mt gA.
$$
The $Sl_n({\mbb C})$-semi-stable points of for this action is 
the set ${\mbb P}(\Hom({\mbb C}^m,{\mbb C}^n))^\sst$ of homomorphisms 
whose rank is $n$. The stabilizer of the stable points is 
${\mbb Z}/n{\mbb Z}$ but this does not represent any problem 
because $PSl_n({\mbb C})$ acts freely, its action linearizes in 
${\cal O}(n)$ and the corresponding semi-stable locus is the 
same as that of the $Sl_n({\mbb C})$-action. The quotient is 
the Grassmannian $Gr_{m-n}({\mbb C}^m)$ with quotient map 
$$
\phi :{({\mbb P}^{mn-1})}^\s\lar Gr_{m-n}({\mbb C}^m),
$$ 
$$
\phi ([A])={\rm Ker}\,A\quad 
\forall A\in {\mbb P}({\Hom}({\mbb C}^m,{\mbb C}^n))^\s.
$$ 
There are morphisms ${\mbb P}^1\rar {({\mbb P}^{mn-1})}^\sst$ in 
each degree $d>0$, an example being
$$
[\zeta _0:\zeta _1]\lmt 
\left(
\begin{array}{cccccc}
\zeta _0^d & \zeta_1^d & \dots & 0 & 0 & 0\\ 
 0 & \ddots & \ddots & \vdots & \vdots & \vdots\\ 
\vdots &\ddots & \ddots & \zeta _1^d & 0 & 0\\ 
 0 & \dots & 0 & \zeta _0^d & \zeta _1^d & 0
\end{array}
\right).
$$
The class induced in the Grassmannian is $nd$ times the class 
of a line. The closure of the inverse image of a point in 
$Gr_{m-n}({\mbb C}^m)$ is a $(n^2-1)$-plane in ${\mbb P}^{mn-1}$; 
this can be seen easily looking at the inverse image of  
$\langle e_{n+1},\dots,e_m\rangle$, where $e_1,\dots,e_m$ 
is the standard basis of ${\mbb C}^m$. Consequently the rational 
slice $\zeta$ for the quotient map is just a $n(m-n)$-plane in the 
projective space, whose class is $H^{n^2-1}$.

The pull-back of a Schubert cycle $\sigma_\lambda$ in the 
Grassmannian is $\phi^*\sigma_\lambda=d(\lambda)H^{|\lambda|}
\in A^{|\lambda|}({\mbb P}^{mn-1})$, for some integer 
$d(\lambda)$.  If $k=\sum_{j=1}^k|\lambda_j|-mnd-n(m-n)+3$, 
we obtain the equality
$$
\begin{array}{l}
GW_{Gr_{m-n}({\mbb C}^m),nd}^k(\sigma_{\lambda_1}
\otimes\dots\otimes\sigma_{\lambda_k})
\\[1ex]\disp 
=\prod_{j=1}^k d(\lambda_j)\cdot 
GW_{{\mbb P}^{mn-1},d}^k(H^{|\lambda_1|+n^2-1}
\otimes\dots\otimes 
H^{|\lambda_k|}).
\end{array}
$$
The question which comes to mind is how can be computed the numbers 
$d(\lambda)$? If $\{0\}\subset F_1\subset\dots\subset F_m={\mbb C}^m$ 
is the standard flag of ${\mbb C}^m$, then 
$$
\phi^*\sigma_\lambda
=\{ [A] \mid \dim({\rm Ker\,}A\cap F_{n+j-\lambda_j})\geq j,\; 
j=1,\dots,m-n\}
$$
and this is just a degeneration locus of the evaluation 
homomorphism 
$$\veps :
{\mbb C}^m\otimes {\cal O}_{{\mbb P}^{mn-1}}(-1)
\lar {\mbb C}^n.
$$
The degree of this subvariety of ${\mbb P}^{mn-1}$ is given by a
determinantal formula which can be found in \cite{F}, theorem 14.3 
page 249. In the particular case when $\sigma_\lambda=\sigma_k$, 
$k=1,\dots,n$, is a special Schubert cycle, we need to compute the 
degree $d_k$ of the subvariety 
$\{ [A] \mid \dim({\rm Ker\,}A\cap F_{n-k+1})\geq 1\}
=\{ [A] \mid \veps_{|F_{n-k+1}}([A])\text{ is not injective}\}
\subset {\mbb P}^{mn-1}$. According to \cite{F}, theorem 14.4 
page 254, $d_k$ is the coefficient of $H^k$ in the development 
of $1/(1-H)^{n-k+1}$, which is $\binom{n}{k}$. 


\section{Second application: Hamiltonian GW-invariants}
{\label{sct:conseq2}}

The second application concerns the so-called Hamiltonian 
Gromov-Witten invariants which were recently introduced in 
\cite{cgs,mun}. The purpose of this and the next sections is 
to put into an algebro-geometric perspective the construction 
performed in these preprints and to show how is that related 
to the problem studied in this article, at least in the case 
of torus actions.

In what follows, $K$ will denote the compact torus ${(S^1)}^r$ 
and $G=K^c={({\mbb C}^*)}^r$ will be its complexification. We 
assume that $G$ acts holomorphically on a projective variety $X$ 
and that the action is linearized in a very ample line bundle 
${\cal O}_X(1)\rar X$. Then the maximal compact subgroup $K$ 
will preserve a symplectic form on $X$ representing the first 
Chern class of ${\cal O}_X(1)$, that is we get a Hamiltonian 
action on $X$. In all the rest, $C$ denotes a smooth projective 
curve with a K\"ahler metric on it. As usual, $\EG\rar\BG$ and 
$\EK\rar\BK$ will stand for the universal $G$ and $K$-bundles; 
they are uniquely determined (up to homotopy) by the condition 
that are contractible and $G$ and $K$ act freely on them, so  
that we may take $\EG=\EK=:E$. 

The $K$-equivariant homology of $X$ is defined as 
$H_*^K(X):=H_*(E\times_K X)$ and elements of it can be constructed 
as follows: one starts with a principal $K$-bundle $P\rar M$ over a 
closed $\crl C^\infty$-manifold $M$ of real dimension $d$ together 
with a $K$-equivariant map $U:P\rar X$. The $d$-dimensional 
equivariant homology class defined by this data is the image of the 
fundamental class of $M$ under 
$$
H_d(M)\srel{\cong}{\llar} H_d^K(P)\lar H_d^K(X).
$$
For every $K$-equivariant $2$-homology class $B\in H_2^K(X;{\mbb Z})$ 
there is a closed Riemann surface $\Sigma$ and a principal $K$-bundle 
$P\rar\Sigma$ together with a $K$-equivariant map $U:P\rar X$ 
representing the class $B$. Moreover, if $\Sigma $ is connected and 
$P,P'\rar\Sigma$ represent the same class, then $P$ and $P'$ are 
isomorphic as $K$-bundles. In other words, the choice of an 
equivariant homology class uniquely determines the topological 
type of the principal bundles, over a {\em fixed} Riemann surface, 
which can represent this homology class. 

Given a principal $K$-bundle $P\rar C$, the complexified bundle 
$P\times_K G$ will be denoted $P^c$. The gauge groups of $P$ and 
$P^c$ are respectively 
$$
\crl G(P)=\{ f:C\rar {(S^1)}^r \}
\quad\text{and}\quad
\crl{G}^c(P)=\{ f:C\rar {({\mbb C}^*)}^r\}.
$$ 
A base point $\zeta_0\in C$ will be fixed once for all. Corresponding 
to it we will consider the based gauges $\crl G_0^c(P)$ of $P^c$ which 
are the identity at $\zeta_0$. The full gauge group is then the direct 
product of $ {({\mbb C}^*)}^r$ with the based gauge group.


\subsection{\sl $A$-holomorphic maps}
{\label{s-sct:A-holo}}

Now we turn to another ingredient used in the definition of the 
Hamiltonian GW-invariants. 

\begin{s-definition}{\label{def:A-holo}}
{\rm (i)} Given a connection $A\in\crl{A}(P)$ and a $K$-equivariant 
map $U:P\rar X$, the operator $\rd_A U$ is defined as
$$
T_pP\ni w\lmt \rd U_p(w)+{\xi(A(w))}_{U(p)},
$$
where ${\xi(a)}_x$ is the tangent vector at $x\in X$ determined by 
$a\in {\it Lie}(K)$.

(ii) A $K$-equivariant map $U:(P,A)\rar X$ is called 
{\em A-holomorphic} if $\dbar_A U=0$, where 
$$
\dbar_A U:=\frac{1}{2}(\rd_A U+J_X\circ\rd _A U\circ J_C).
$$
The notation $J_C$ stands for the complex structure induced on the 
$A$-horizontal spaces of $P$ by the complex structure of $C$. 
\end{s-definition}

In more down to earth terms, a $K$-equivariant map 
$U:(P,A)\rar X$ is $A$-holomorphic if and only if 
$$
\rd U_p(\wtld{Jv})=J\rd U_p(\wtld{v}),
$$
where $\wtld{v}$ and $\wtld{Jv}$ denote respectively 
the $A$-horizontal liftings in $p\in P$ of the vectors 
$v, Jv\in TC$. 

\begin{s-definition}(cf. \cite{cgs} section 3.2) 
Denote 
$$
{\mfrak X}_{B}:=\left\{ (U,A)\in C_K^\infty(P,X;B)\times 
\crl{A}(P) \mid \dbar_A U=0\right\}
$$
the space of $K$-equivariant, $A$-holomorphic smooth maps which 
represent the class $B\in H_2^K(X)$. 
\end{s-definition}

Any $K$-equivariant map $U:P\rar X$ induces a map 
$\bar u:C\rar P\times_K X$. Then $U$ is $A$-holomorphic if and 
only if $\bar u$ is holomorphic. One has to be careful with the 
(integrable) complex structure on $P\times_K X$ which is induced by 
the connection $A$. For vectors tangent to the fibres of 
$P\times_K X\rar C$ the complex structure agrees with that of $X$, 
while for $v\in TC$ (in a local trivialization $P\cong C\times K)$,
\begin{align}{\label{eq:A-cplx-struct}}
J_{P\times_K X}(v)=Jv+\xi (A(Jv))-J\xi(A(v)).
\end{align}
In this formula and in all the rest of the paper, for 
$a\in{\it Lie}(K)$, $\xi(a)$ will denote the vector field 
on $X$ induced by the infinitesimal $K$-action.

Clearly, the real gauges act on ${\mfrak X}_{B}$ but it turns 
out that the complex gauges act also. The formula for it is given 
in the lemma below. 

\begin{s-lemma}{\label{lm:gauge-action}}
The complex gauges $\crl{G}^c(P)$ act on ${\mfrak X}_{B}$ as  
$f\times (U,A)\lmt (f\cdot U,f\cdot A)$, where 
\begin{eqnarray}{\label{eq:fU}}
(f\cdot U)(p):=f^{-1}(p)U(p)
\end{eqnarray} 
and 
\begin{eqnarray}{\label{eq:fA}}
f\cdot A:=A+{(f^{-1}\rd f)}_{\mfrak k}
+*{(f^{-1}\rd f)}_{i\mfrak k}.
\end{eqnarray}
\end{s-lemma}

Some explanation is in order: any 
$a\in {\it Lie}(G)={\it Lie}(K)\oplus i{\it Lie}(K)$ 
can be uniquely written 
$a=a_{\mfrak k}+ia_{i\mfrak k}$, 
with 
$a_{\mfrak k},a_{i\mfrak k}\in {\it Lie}(K)$. 
The $*$ in the formula represents the Hodge star operator on $C$.

\begin{proof}
It is clear that formula \eqref{eq:fU} 
just extends the action of the real gauges by composition 
on the right. We shall prove the formula for the action of $f$ 
on $A$ searching a connection $A'$ on $P$ which makes 
$A'$-holomorphic the map $U':=f\cdot U$. For doing computations 
we use a local trivialization of $P$, so that $P$ itself may be 
assumed trivial (as long as the objects found in the end are 
globally defined).

In what follows, $\zeta$ denotes a point on $C$. By assumption 
$P=C\times K$ and $U(\zeta ,g^{-1})=gu(\zeta)$, for some 
$u:C\rar X$. I want to find a connection $A'$ on $P$ such that 
$\dbar_{A'} U'=0$. Since $U'$ is $K$-equivariant, it is enough to 
check this condition at points $(\zeta,1)\in P$. Because $U$ is 
$A$-holomorphic, 
$$
\rd U_{(\zeta,1)}(Jv-A(Jv))=J\rd U_{(\zeta ,1)}(v-A(v))
$$
for $v\in T_\zeta C$, or equivalently 
\begin{align}{\label{eq:A-hol}}
\rd u_\zeta (Jv)+{\xi(A(Jv))}_{u(\zeta)}=
J\rd u_\zeta (v)+J\xi (A(v))_{u(\zeta)}.
\end{align} 
Formula \eqref{eq:fU} implies that 
$$
\begin{array}{rl}
\rd U'_{(\zeta ,1)}(Jv-A'(Jv))&
=f^{-1}(\zeta)\rd u_\zeta (Jv)
-f^{-1}(\zeta){\xi((f^{-1}\rd f)(Jv))}_{u(\zeta)}\\ 
&
+f^{-1}(\zeta){\xi(A'(Jv))}_{u(\zeta)}
\end{array}
$$ 
and 
$$
\begin{array}{rl}
J\rd U'_{(\zeta ,1)}(v-A'(v))&
=f^{-1}(\zeta)J\rd u_\zeta (v)
-f^{-1}(\zeta)J{\xi ((f^{-1}\rd f)(v))}_{u(\zeta)}\\ 
&
+f^{-1}(\zeta)J{\xi(A'(v))}_{u(\zeta)}.
\end{array}
$$ 
For $U'$ to be $A'$-holomorphic it is necessary and sufficient 
that the difference of these two quantities is zero. Imposing 
this condition, we find 
$$
\begin{array}{rl}
0&
=J\rd u_\zeta(v)-\rd u_\zeta (Jv)
-J{\xi((f^{-1}\rd f)(v))}_{u(\zeta)}
+{\xi((f^{-1}\rd f)(Jv))}_{u(\zeta)}\\ 
&
+J{\xi(A'(v))}_{u(\zeta)}-{\xi(A'(Jv))}_{u(\zeta)}\\[1.5ex] 
&  
\srel{_{\eqref{eq:A-hol}}}{=}
{\xi(A(Jv))}_{u(\zeta)} -J{\xi (A(v))}_{u(\zeta)}
-{\xi(A'(Jv))}_{u(\zeta)}
+J{\xi(A'(v))}_{u(\zeta)}\\ 
&
+{\xi((f^{-1}\rd f)(Jv))}_{u(\zeta)}
-J{\xi((f^{-1}\rd f)(v))}_{u(\zeta)}\\[1.5ex] 
&  
={\xi\bigl( A(Jv)-A'(Jv)\bigr) }_{u(\zeta)}
-J{\xi\bigl( A(v)-A'(v)\bigr) }_{u(\zeta)}\\
&
+{\xi\bigl( (f^{-1}\rd f)(Jv)\bigr) }_{u(\zeta)}
-J{\xi\bigl( (f^{-1}\rd f)(v)\bigr) }_{u(\zeta)}.
\end{array}
$$
It remains to separate the ${\it Lie}(K)$ and $i{\it Lie}(K)$ 
components of the last line. 
$$
\begin{array}{l}
\xi\bigl( (f^{-1}\rd f)(Jv)\bigr)
-J\xi\bigl( (f^{-1}\rd f)(v)\bigr)
\\[1ex]
=\xi\bigl( (f^{-1}\rd f)_{\mfrak k}(Jv)
+i(f^{-1}\rd f)_{i\mfrak k}(Jv)\bigr)
-\xi\bigl( i(f^{-1}\rd f)_{\mfrak k}(v)
-(f^{-1}\rd f)_{i\mfrak k}(v)\bigr)
\\[1ex]
=\xi
\bigl( (f^{-1}\rd f)_{\mfrak k}(Jv)
+(f^{-1}\rd f)_{i\mfrak k}(v)\bigr)
-J\xi\bigl( (f^{-1}\rd f)_{\mfrak k}(v)
-(f^{-1}\rd f)_{i\mfrak k}(Jv)\bigr).
\end{array}
$$
Inserting this into the previous relation, we obtain
$$
\begin{array}{rl}
0\kern-1ex&
=\xi\bigl( A(Jv)-A'(Jv)+(f^{-1}\rd f)_{\mfrak k}(Jv)
+(f^{-1}\rd f)_{i\mfrak k}(v)\bigr)\\ 
&
-J\xi\bigl( A(v)-A'(v)+(f^{-1}\rd f)_{\mfrak k}(v)
-(f^{-1}\rd f)_{i\mfrak k}(Jv)\bigr).
\end{array}
$$
For $A'$ defined by 
$$
A'=A+(f^{-1}\rd f)_{\mfrak k}-(f^{-1}\rd f)_{i\mfrak k}\circ J,
$$
the last equality is satisfied. Notice that in general this is the 
only possible choice for $A'$ since the vectors $\xi $ and $J\xi$ 
are linearly independent in most cases. 

Using local normal coordinates on $C$, it follows that for any 
$1$-form $\alpha\in \Omega_C^1$, $\alpha\circ J=-(*\alpha)$.
\end{proof}

\begin{s-remark}{\label{rk:f=exp}} 
(i) It follows from  formula \eqref{eq:fA} that the 
$\crl{G}_0^c(P)$-action on ${\mfrak X}_{B}$ is free. 

(ii) Any $f:C\rar {({\mbb C}^*)}^r$ is of the form 
$f(\zeta)=R(\zeta)\vphi(\zeta)$, with $R:C\rar {\mbb R}^r$ 
and $\vphi:C\rar {(S^1)}^r$. Formula \eqref{eq:fA} becomes 
\begin{eqnarray}{\label{eq:fxA}}
f\cdot A=A+\vphi ^{-1}\rd \vphi -i*\rd(\log R).
\end{eqnarray}
The form $\vphi^{-1}\rd\vphi$ is closed, but not necessarily exact; 
it is exact if and only if $\vphi_*:\pi_1(C)\rar {(\pi_1(S^1))}^r$ 
is the zero homomorphism, or equivalently $\vphi =\exp(i\theta)$ 
for some $\theta:C\rar{\mbb R}^r$. However it always defines an 
integral $1$-cohomology class and conversely, any integral 
$1$-cohomology class can be represented in this form. 

Using the Hodge decomposition of $\Omega_C^1$, this discussion 
implies that the pointed complex gauge equivalence classes of 
connections in $P$ are parameterized by
$$
{\mbb H}^1(C,{\mbb R}^r)/{\mbb H}^1(C,{\mbb Z}^r),
$$
where ${\mbb H}^1(C,{\mbb R}^r)$ denotes the space of harmonic 
${\mbb R}^r$-valued $1$-forms on $C$. This quotient is just the 
$r^{\rm th}$ power of the familiar Picard variety of $C$, when 
${\mbb H}^1(C,{\mbb R}^r)$ is given the complex structure defined 
by the Hodge-star of $C$.  

(iii) In the genus zero case, i.e. $C\cong {\mbb P}^1$, all 
gauges admit a globally defined logarithm. Therefore all connections 
are gauge equivalent, which is the same saying that in a given 
topological principal bundle $P\rar{\mbb P}^1$ there is only one 
equivalence class of holomorphic structures.
\end{s-remark}


\subsection{\sl Short digression on the Picard variety}
{\label{s-sct:picard}}

All the statements in this section should be well known, but 
we are recalling them in order to fix the notations. As I have 
already mentioned, the quotient 
$$
{\mbb H}^1(C,{\mbb R}^r)/{\mbb H}^1(C,{\mbb Z}^r)
={\left( {\mbb H}^1(C,{\mbb R})/{\mbb H}^1(C,{\mbb Z})\right) }^r
={(Pic^0C)}^r
$$
is the $r^{\rm th}$ power of the Picard variety of $C$, when 
${\mbb H}^1(C,{\mbb R}^r)$ is regarded as a complex vector space 
with complex structure given by the Hodge-star of $C$. It is a 
projective torus which parameterizes topologically trivial, 
holomorphic principal ${({\mbb C}^*)}^r$-bundles over $C$. 

Let 
$$
{\cal P}_r^c\lar {(Pic^0C)}^r\times C
$$
be a universal principal ${({\mbb C}^*)}^r$-bundle. It has the 
property that for any point $\tau\in {(Pic^0C)}^r$, the restriction 
${({\cal P}_r^c)}_\tau\lar C$ represents the point $\tau$. We shall 
be interested in describing a connection in this bundle which induces 
its holomorphic structure. 

Let us start with the 

\begin{s-lemma}{\label{lm:easy}}
For any $\alpha\in {\mbb H}^1(C;{\mbb Z}^r)$, there is a unique 
$\vphi_\alpha:C\rar {(S^1)}^r$ such that $\vphi_\alpha(\zeta_0)=1$ 
and $\vphi_\alpha^{-1}\rd\vphi_\alpha=\alpha$. (The point 
$\zeta_0\in C$ was fixed from the beginning). 
\end{s-lemma}

\begin{proof} Clearly, we may assume that $r=1$. The uniqueness part 
is immediate. For the existence part, notice that if $\alpha$ is 
exact, i.e. $\alpha=\rd\theta$ for $\theta :C\rar {\mbb R}$, then 
$\vphi:=\exp\left(i(\theta-\theta(\zeta_0))\right)$ does the job. 
Homotopy classes of maps $C\rar{\mbb R}$ are parameterized by 
${\rm Hom}_{\mbb Z}(H^1(C;{\mbb Z});{\mbb Z})\cong H^1(C;{\mbb Z})$, 
so for $\alpha\in {\mbb H}^1(C;{\mbb Z})$ there exists 
$\vphi':C\rar S^1$ such $\vphi'(\zeta _0)=1$ and 
$[{(\vphi ')}^{-1}\rd\vphi']=[\alpha]\in H^1(C;{\mbb Z})$. By the 
discussion above, there exists $\vphi_0:C\rar S^1$ such that 
$\vphi _0(\zeta _0)=1$ and $\vphi _0^{-1}\rd\vphi_0=
\alpha-{(\vphi ')}^{-1}\rd\vphi'$. Now $\vphi:=\vphi_0\vphi '$ will be 
convenient. 
\end{proof}

\begin{s-remark}{\label{rk:morph}} 
The map 
$$
{\mbb H}^1(C;{\mbb Z}^r)\ni\alpha\lmt \vphi _\alpha\in 
{\crl C}^\infty\bigl( C,{(S^1)}^r\bigr)
$$
is a morphism of groups i.e. $\vphi_\alpha\vphi_\beta
=\vphi_{\alpha+\beta}$.
\end{s-remark}

Fix once for all a real connection $A_0$ in the smooth 
$({\mbb C}^*)^r$-bundle $P^c\rar C$ i.e. one coming from 
a connection in the real ${(S^1)}^r$-bundle. 

\begin{s-lemma}{\label{lm:holo-str}}
(i) On ${\mbb H}^1(C;{\mbb R}^r)\times C$, there is a natural, 
closed ${\mbb R}^r$-valued $1$-form $\chi$ which is defined by
$$
\chi_{(A,\zeta )}(a,v):=A(v)\quad{\rm for}\quad(a,v)\in 
T_{(A,\zeta )}\bigl({\mbb H}^1(C;{\mbb R}^r)\times C\bigr).
$$
(ii) The (real) connection ${{\cal A}}:=A_0+\chi$ defines a 
holomorphic structure on the bundle
${\rm pr}_C^*P^c\rar {\mbb H}^1(C;{\mbb R}^r)\times C$.
\end{s-lemma}

\begin{proof} The curvature of ${\cal A}$ is $F_{{\cal A}}
={\rm pr}_C^*F_{A_0}+\rd\chi$. For $a\in\mbb{H}^1(C;\mbb{R}^r)$ 
and $v\in TC$, $\rd\chi\bigl((a,0),(0,v)\bigr)=a(v)$, 
and $\rd\chi$ evaluates zero on other pairs of vectors. 

It is easy to see that for any $1$-form $a$ on $C$ and any 
tangent vector $v$ to $C$, $(*a)(Jv)=a(v)$. This implies that 
the connection ${\cal A}$ defines indeed a holomorphic structure 
in ${\rm pr}_C^*P^c$ because its curvature is a $(1,1)$-form 
on ${\mbb H}^1(C;{\mbb R}^r)\times C$.
\end{proof}

\begin{s-proposition}{\label{prop:univ-bdl}}
(i) The group ${\mbb H}^1(C;{\mbb Z}^r)$ acts holomorphically, 
by real gauges, on ${\rm pr}_C^*P^c$ by
$$
{\mbb H}^1(C;{\mbb Z}^r)\times 
\bigl( {\mbb H}^1(C;{\mbb R}^r)\times P^c\bigr)
\lar {\mbb H}^1(C;{\mbb R}^r)\times P^c,
$$ 
$$
\alpha\times (A,p):=(A+\alpha, R_{\vphi _\alpha}p).
$$

(ii) The holomorphic principal bundle 
$$
{\cal P}_r^c(A_0):={\rm pr}_C^*P^c/{\mbb H}^1(C;{\mbb Z}^r)
\lar {(Pic^0C)}^r\times C
$$
is a universal principal ${({\mbb C}^*)}^r$-bundle which 
parameterizes holomorphic bundles over $C$ having fixed topological 
type defined by $A_0$. It also comes with the connection induced by 
${\cal A}$. 
\end{s-proposition}

\begin{proof} (i) Remark \ref{rk:morph} implies that the formula 
above is indeed an action. It is also holomorphic because 
${\mbb H}^1(C;{\mbb Z}^r)$ preserves the connection ${\cal A}$; indeed, 
$A+\alpha=\vphi_\alpha A.$

(ii) The statement is a direct consequence of \ref{rk:f=exp}. 
\end{proof}

 We should say that we have worked with {\em complex} principal 
bundles throughout this section because the accent was put on 
their {\em holomorphic} structure. But ${\cal P}_r^c(A_0)$ is the 
complexification of a real ${(S^1)}^r$-bundle ${\cal P}_r(A_0)\rar 
{(Pic^0C)}^r\times C$ and the connection ${\cal A}$ comes from a 
connection in ${\cal P}_r(A_0)$, because the connection $A_0$, 
fixed from the beginning, was {\em real} and the action on 
${\rm pr}_C^*P^c$ was done by {\em real} gauges.


\section{Moduli spaces}
{\label{sct:spaces}}

In this section we shall see that the spaces introduced in \cite{cgs} 
and \cite{mun} for defining invariants of Hamiltonian group actions 
have nice algebraic interpretation.

The authors of \cite{cgs} introduce (a perturbation of) 
\begin{align}{\label{tlds-space}}
\tld S_{C,k}(X;B):=\left.
\bigl\{ (U,A)\in {\mfrak X}_{B} \mid *F_A+m\circ U=0\bigr\}
\times {(P^k)}_o\right/ \crl{G}(P)
\end{align}
and
\begin{align}{\label{s-space}}
S_{C,k}(X;B):=\tld S_{C,k}(X;B)/G^k.
\end{align}
In the definition above $m:X\rar {\mbb R}^r$ is the moment map 
corresponding to the ${(S^1)}^r$-action and ${(P^k)}_o$ is the 
open set in $P^k$ consisting of $k$ points which project to 
$k$ distinct points of $C$. The group actions are as 
follows:
$$ 
f\times \bigl( (U,A)\times(p_1,\dots,p_k)\bigr) :=
(f\cdot U,f\cdot A)\times(R_{f(p_1)}p_1,\dots,R_{f(p_k)}p_k)
$$
and
$$
(g_1,\dots,g_k)\times 
\bigl[ (U,A)\times(p_1,\dots,p_k)\bigr] :=
\bigl[ (U,A)\times(R_{g_1}p_1,\dots,R_{g_k}p_k)\bigr],
$$
for $f\in\crl{G}(P)$ and $g_1,\dots,g_k\in K={(S^1)}^r$. The letter 
`$R$' denotes the right action of $G$ on the principal bundle. The 
expected dimension of this space is 
\begin{align}{\label{eq:exp-dim}}
2D:={\rm exp.dim}_{\mbb R}S_{C,k}(X;B)=
2(1-g)(n-r)+2c_1^K(X)\cdot B+2k,
\end{align}
where $c_1^K(X)$ denote the $K$-equivariant first Chern class 
of $X$.

The space $S_{C,k}(X,A)$ can also be expressed as the infinite 
dimensional invariant quotient as 
\begin{eqnarray}{\label{git-s-space}}
S_{C,k}(X,B)=
\left.
\left({\mfrak X}_B^{\rm s}\times {(P^c)}_o^k/\crl{G}^c(P)\right)
\right/ G^k.
\end{eqnarray}
The notation ${\mfrak X}_B^{\rm s}\subset{\mfrak X}_{B}$ stands for 
the set of so-called stable pairs and $\crl{G}^c(P)$ acts on 
${\mfrak X}_{B}$ as described in lemma \ref{lm:gauge-action}.

For the present purposes, it will be more convenient to use this 
second description. Since $\crl{G}^c(P)=\crl{G}_0^c(P)\times G$ 
and the actions of $\crl{G}^c(P)$ and $G^k$ on 
${\mfrak X}_{B}$ commute, $S_{C,k}(X;B)$ can be constructed in a 
different way: first take the quotient for the free 
$\crl{G}_0^c(P)\times G^k$-action (which is finite dimensional) 
and after make the invariant quotient for the remaining $G$-action. 

$K$-equivariant maps $U:P\rar X$ and $G$-equivariant maps 
$U^c:P^c\rar X$ will be used interchangeably. That there is no 
harm in doing so follows from the fact that any $K$-equivariant 
$U$ defines (in the obvious way) a $G$-equivariant $U^c$; 
conversely, any such $U^c$ defines a corresponding $U$ composing 
it with the inclusion $j:P\hra P^c$. 

\begin{lemma}{\label{lm:u-bar}}
(i)  The variety $\bar X$ defined by 
$$
\bar X:={\cal P}_r^c(A_0)\times_G X={\cal P}_r(A_0)\times_K X, 
$$
carries a natural structure of a complex projective variety. Its 
complex dimension is $\dim \bar X=gr+\dim X+1$, where $g$ is the 
genus of $C$ and $r=\dim G$. 

(ii) Any $K$-equivariant, $A$-holomorphic map $U:P\rar X$, with 
$A\in \crl{A}(P)$, which represents an equivariant $2$-homology 
class $B\in H_2^K(X;{\mbb Z})$ defines a holomorphic map 
$\bar u:C\rar\bar X$ which represents a class 
$\bar B\in H_2(\bar X;{\mbb Z})$ depending on $B$ only. 
If $\pi:\bar X\rar{(Pic^0C)}^r\times C$ is the natural projection, 
$\pi_*\bar B=0\oplus [C]$. Moreover, 
$c_1^K(X)\cdot B=c_1(T_{\bar X}^{\rm rel})\cdot \bar B$, where 
$T_{\bar X}^{\rm rel}$ denotes the $\pi$-relative tangent bundle.

(iii) Consider an $A$-holomorphic, $K$-equivariant map 
$U:P\rar X$ and $g\in\crl{G}_0^c(P)$. Then $U$ and $gU$ 
define the same map $C\rar\bar X$. 
\end{lemma}

\begin{proof}
(i) $\bar X$ has a holomorphic structure because 
${\cal P}_r^c(A_0)\rar {(Pic ^0C)}^r\times C$ is a holomorphic 
bundle, according to lemma \ref{lm:holo-str}. That it is also a 
projective variety follows from the fact that the Picard torus 
is projective. 

(ii) Remark \ref{rk:f=exp} implies that given $U:(P,A)\rar X$ 
there is a unique $f_A\in\crl{G}_0^c(P)$ (depending on $A$) such 
that 
$$
f_A\cdot A=\text{harmonic part of }A-A_0=:h(A-A_0)\in 
{\mbb H}^1(C;{\mbb R}^r).
$$ 
The composed map
\begin{align}{\label{comp-map}}
P
\xrar{(\{ h(A-A_0)\}\times{\rm id}_P)\times f_AU} 
\underbrace{{\mbb H}^1(C;{\mbb R}^r)\times P}_{={\rm pr}_C^*P}
\times X
\xrar{\phantom{MMM}} 
{\cal P}_r(A_0)\times X
\end{align}
is $K$-equivariant and therefore defines
$$
\bar u:C\lar {\cal P}_r(A_0)\times_K X=\bar X.
$$

Since $f_A U:P\rar X$ is $f_A A$-holomorphic, it follows 
that this map is holomorphic. If $p_\zeta\in P$ (or in $P^c$) 
denotes a point lying over $\zeta\in C$, the explicit formula for 
$\bar u$ is
\begin{eqnarray}{\label{eq:u-bar}}
\zeta\srel{\bar u}{\lmt }
\bigl[ [h(A-A_0),p_\zeta], (f_AU)(p_\zeta)\bigr],
\end{eqnarray}
where the square brackets denote obvious equivalence classes.

Suppose that $P=\rho^*E$ for a map $\rho:C\rar\BK$. Then from 
the commutative diagram
$$
\xymatrix{
T_{P\times_K X}^{\rm rel}=\kern-2.5em & 
{(\rho ,{\rm id}_X)}^*(E\times_K T_X)
\ar[d]\ar[r] &
E\times_K T_X \ar[d]\\ 
C\ar[r]^-{\bar u} & 
P\times_K X \ar[r]^-{(\rho,{\rm id}_X)}&
E\times_K X.
}
$$
we can see that the homology class $\bar B$ depends only on the 
$K$-equivariant class $B$ and also that 
$c_1^K(X)\cdot B=c_1(T_{\bar X}^{\rm rel})\cdot \bar B$. 

(iii) Consider now $U:P\rar X$ and $g\in \crl{G}_0^c(P)$, 
$g=R\cdot\vphi$ with $\vphi:C\rar {(S^1)}^r$ and 
$R:C\rar {\mbb R}^r$. According to \eqref{eq:fxA}, 
$gA=A+\vphi^{-1}\rd\vphi-i*\rd (\log R)$, so 
$$
h(gA-A_0)=h(A-A_0)+h(\vphi^{-1}\rd\vphi).
$$
Notice that $\alpha _\vphi:=h(\vphi^{-1}\rd\vphi)$ is actually an 
{\em integral} ${\mbb R}^r$-valued harmonic form and according to 
lemma \ref{lm:easy} there is a unique $\psi_\vphi\in \crl{G}_0(P)$ 
such that $\psi_\vphi^{-1}\rd\psi_\vphi=\alpha_\vphi$. I claim that 
$f_{gA}=\psi_\vphi f_A g^{-1}$, {\it i.e.} that 
$(\psi_\vphi f_A g^{-1})(gA)=A_0+h(gA-A_0)$. Indeed,
$$
\begin{array}{rl}
(\psi_\vphi f_A g^{-1})(gA)&=\psi_\vphi (f_A A)
=\psi_\vphi (A_0+h(A-A_0))\\ 
&
=A_0+h(A-A_0)+\psi_\vphi^{-1}\rd\psi_\vphi\\ 
&
=A_0+h(A-A_0)+\alpha _\vphi=A_0+h(gA-A_0).
\end{array}
$$
We are going to check now that $U$ and $gU$ induce the same 
(holomorphic) map $C\rar\bar X$. Using formula \eqref{eq:u-bar},  
$gU$ reads 
\begin{align}{\label{well-def}}
\begin{array}{crl}
\zeta& 
\lmt& 
\bigl[ [h(A-A_0)+\alpha_\vphi ,p_\zeta ],
\psi_\vphi^{-1}(p_\zeta)\cdot(f_A U)(p_\zeta)\bigr]\\[1ex] 
&=
&\bigl[ [h(A-A_0)+\alpha_\vphi ,R_{\psi_\vphi}p_\zeta],
(f_A U)(p_\zeta)\bigr]\\[1ex] 
&= 
&\bigl[ [h(A-A_0),p_\zeta],(f_A U)(p_\zeta)\bigr].
\end{array}
\end{align}
This finishes the proof of the lemma.
\end{proof}

The next proposition gives the algebro-geometric interpretation of 
the space \eqref{s-space}. 

\begin{proposition}{\label{prop:1-1}}
There is a one-to-one 
map 
$$
{\mfrak X}_{B}\times {(P^c)}_o^k/\crl{G}_0^c(P)\times{G}^k 
\leftarrow\kern-1em\xrar{\phantom{x}1:1\phantom{x}}
M_{C,k}(\bar X;\bar B),
$$
where $M_{C,k}(\bar X;\bar B)$ denotes the space of stable maps 
$(C,\unl x)\rar\bar X$ with $k$ marked points $\unl x\in{(C^k)}_o$ 
and representing the $2$-homology class $\bar B$. As usual, 
${(C^k)}_o$ denotes the complement in $C^k$ of the diagonals.
\end{proposition}

\begin{proof} Consider an $A$-holomorphic, $K$-equivariant map 
$U:P\rar X$ together with $k$ marked points 
$\unl p\in{(P^c)}^k_o$ and let 
$\unl x:=\pi_{_C}(\unl p)\in{(C^k)}_o$. 
Lemma \ref{lm:u-bar} says that this data induces a morphism 
$C\rar\bar X$ and moreover, it does not depend on the 
$\crl{G}_0^c(P)$-orbit of $(U,A)$. So we get a map
$$
F:{\mfrak X}_{B}\times {(P^c)}_o^k/\crl{G}_0^c(P)\lar 
M_{C,k}(\bar X;\bar B)
$$ 
This map is clearly ${G}^k$-invariant and therefore descends 
to the quotient 
$$
F:{\mfrak X}_{B}\times {(P^c)}_o^k/\crl{G}_0^c(P)\times G^k 
\lar M_{C,k}(\bar X;\bar B).
$$
Because the composition $\pi\bar u=\{\tau\}\times{\rm id}_C$, 
for some $\tau\in {(Pic^0C)}^r$, the map $\bar u$ is in fact a 
representative for the corresponding stable map (see definition 
\ref{def:stable-map}). 

The map $F$ is clearly surjective: given a point $(\bar u,\unl x)
\in M_{C,k}(\bar X;\bar B)$, consider the diagram
$$
\xymatrix{
{\bar u}^*({\cal P}_r(A_0)\times X)=\kern-2.6em & 
P\ar[d] \ar[r]^-{\bar U} & 
{\cal P}_r(A_0)\times X 
\ar[d] \ar[r]^-{{\rm pr}_X} & X\\ 
 & C\ar[r]^-{\bar u} & \bar X & 
}
$$
The composed map $U:={\rm pr}_X\circ \bar U$ will be a 
$K$-equivariant, $A$-holomorphic map, for $A:={\bar u}^*{{\cal A}}$ 
(see \ref{lm:holo-str} for the definition of ${\cal A}$). As marked 
points in $P$, one may take any $\unl p$ lying over $\unl x$.  

We have to prove that $F$ is injective. Consider  
$U:(P,A,\unl p)\rar X$ and $U':(P,A',\unl p')\rar X$ which induce 
the same morphism $\bar u:C\rar \bar X$. It follows from definition 
\eqref{comp-map} that necessarily 
$h(f_A A-A_0)\equiv h(f_{A'}A'-A_0)\ {\rm mod}\ 
{\mbb H}^1(C;{\mbb Z}^r)$, so that 
$h(f_A A-f_{A'}A')\equiv 0\ {\rm mod}\ {\mbb H}^1(C;{\mbb Z}^r)$. 
Remark \ref{rk:f=exp} implies that $f_A A$ and $f_{A'}A'$ are 
in the same $\crl{G}_0^c(P)$-orbit and consequently $A$ and $A'$ 
are also in the same $\crl{G}_0^c(P)$-orbit. Since $\bar u$ is 
gauge invariant, we may assume that $A=A'$ and even that $A-A_0$ 
is a harmonic form. 

The problem is reduced to the following: two maps 
$U:(P,A,\unl p)\rar X$ and $U':(P,A,\unl p')\rar X$ which define 
the same $(\bar u,\bar x)$ must be equal. Formula \eqref{eq:u-bar} 
says that
$$
\bigl[ [A-A_0,p_\zeta],U(p_\zeta),[\unl p]\bigr]=
\bigl[ [A-A_0,p_\zeta],U'(p_\zeta),[\unl p']\bigr] 
\quad \forall p_\zeta\in P.
$$ 
A moment's thought shows that this imply $U=U'$ and 
$[\unl p]=[\unl p']$.
\end{proof}

\begin{remark}{\label{rk:probl}}
An advantage of working with $M_{C,k}(\bar X;\bar B)$ is that 
it has a natural quasi-projective scheme structure. This was 
proved in lemma \ref{lm:quasi-proj-global-cover}, where is 
also described the construction of its compactification 
$\ovl M_{C,k}(\bar X;\bar B)$ in terms of stable maps. 
Certainly, working within this algebraic frame has its own 
disadvantages: the space of stable maps may be badly behaved 
or it may have wrong dimension; an instance of a very unpleasant 
situation is when 
$M_{C,k}(\bar X;\bar B)\subset\ovl M_{C,k}(\bar X;\bar B)$ 
is not dense. 

When $\ovl M_{C,k}(\bar X;\bar B)$ has larger dimension than 
the expected one, it seems possible to introduce a virtual 
class on $\ovl M_{C,k}(\bar X;\bar B)$ using obstruction 
theory {\em relative to} $\pi:\bar X\rar{(Pic^0C)}^r\times C$. 
This should correspond to the limit of the fundamental cycles 
of the moduli spaces of pseudo-holomorphic curves to 
$(\bar X,\bar J):={\cal P}_r(A_0)\times_K (X,J)$, with $J$ a 
$K$-invariant, generic almost complex structure on $X$. 
Unfortunately, for the moment, I can not make this statement 
more precise.
\end{remark}

There is a natural $G$-action on $\bar X$: 
$$
g\times \bigl[ [A,p],x\bigr]:=\bigl[ [A,p],gx\bigr]
\quad{\rm for}\quad
\bigl[ [A,p],x\bigr] \in \bar X={\cal P}_r^c(A_0)\times_G X
$$ 
which is well-defined precisely because $G$ is commutative.
The $G$-action can be linearized in the line bundle 
$$
L:={\cal P}_r^c(A_0)\times_{G} {\cal O}_X(1)\lar \bar X,
$$ 
which is $\pi$-ample. For $\ell \rar {(Pic^0C)}^r\times C$ 
sufficiently ample line bundle, 
$$
\bar L:=\pi^*\ell\otimes L\lar \bar X
$$
is ample and the $G$-action linearizes again. In is rather 
clear that the set of $G$-semi-stable points of $\bar X$ for 
this action is $\bar X^\sst={\cal P}_r^c(A_0)\times_{G} X^\sst$. 
In particular, if $G$ acts freely on $X^\sst$, it will does the 
same on $\bar X^\sst$. The invariant quotient is, in any case, 
$$
\bar X\invq G={(Pic^0 C)}^r\times C\times \hat X
\quad{\rm where}\quad\hat X:=X\invq G.
$$ 

The next lemma is useful to ``visualize'' better 
$M_{C,k}(\bar X;\bar B)$.

\begin{lemma}{\label{lm:tau}}
Assume that $G$ acts freely on $X^\sst$ and let 
$\phi:X^\sst\rar\hat X$ be the quotient map. Consider 
$\bar u\in M_{C,k}(\bar X;\bar B)$ with 
${\rm Image}(\bar u)\subset{{\cal P}_r^c(A_0)}_\tau\times_GX^\sst$. 
Then ${(\phi\circ\bar u)}^*X^\sst\rar C$ represents the point 
$\tau\in{(Pic^0C)}^r$. 
\end{lemma}

\begin{proof} Notice that in the diagram 
$$
\xymatrix{
 & {{\cal P}_r^c(A_0)}_\tau\times X^\sst
\ar[d] &  X^\sst\ar[d]\\ 
C\ar[r]^-{\bar u} 
& 
{{\cal P}_r^c(A_0)}_\tau\times_GX^\sst
\ar[r]^-{\phi} 
& \hat X,
}
$$
${{\cal P}_r^c(A_0)}_\tau\times X^\sst=\phi^*X^\sst$. Indeed, 
for $\bigl( [p,x],x'\bigr)\in \phi^*X^\sst$ there is a unique 
$p'$ such that $[p',x']=[p,x]$, so we may identify 
$\bigl( [p,x],x'\bigr)=(p',x')$. 
Consequently, ${(\phi\circ\bar u)}^*X^\sst
=\bar u^*({{\cal P}_r^c(A_0)}_\tau\times X^\sst)$ 
and we obtain a $G$-equivariant map 
${(\phi\circ\bar u)}^*X^\sst\rar {{\cal P}_r^c(A_0)}_\tau$ 
which covers the identity of $C$. This one must be an isomorphism.
\end{proof}

We should recall that for obtaining the moduli space 
$S_{C,k}(X;B)$ we still need to divide out the remaining 
$G$-action on $\ovl M_{C,k}(\bar X;\bar B)$. For comparing the 
two points of view, the real-analytic and the algebraic one, we 
shall use the results obtained in section \ref{sct:symp-pers}: 
the moment map corresponding to the $G$-action on $\bar X$ is 
the function $\bar m$ defined by 
$$
\xymatrix{
{\cal P}_r^c(A_0)\times X
\ar[d]\ar[r]^-{{\rm pr}_X}
& X\ar[r]^-{m}
& {\mbb R}^r\\ 
\bar X\ar[urr]_-{\bar m} & 
}
$$
It follows from remark \ref{rk:limit-mm} that the limit 
moment map on $\ovl M_{C,k}(\bar X;\bar B)$ is 
$$
\bar u\lmt \int_C \bar m\circ\bar u\in {\mbb R}^r
$$
and the invariant quotient is constructed dividing the zero 
level set by $K={(S^1)}^r$. 

On the other hand, the zero level $\{*F_A+m\circ U=0\}$ 
(modulo gauge) appearing in the definition \ref{tlds-space} 
can be written also as 
$$
U\lmt \int_C m\circ U=-\int_C F_A=\unl\delta\in{\mbb R}^r,
$$ 
where $\unl\delta$ represents the multi-degree of the 
${(S^1)}^r$-bundle $P$, which is a topological invariant. 
We deduce that $S_{C,k}(X;B)$ defined by \eqref{s-space} 
and $\ovl M_{C,k}(\bar X;\bar B)\invq G$ should be 
birational because both of them are Marsden-Weinstein 
quotients of $\ovl M_{C,k}(\bar X;\bar B)$. However, unless 
$\ovl M_{C,k}(\bar X;\bar B)$ is irreducible, this issue can 
be quite tricky. 

In the particular case when $\unl\delta=0$, that is we are 
looking at topologically trivial bundles, the real-analytic 
and the algebro-geometric constructions coincide.


\section{Hamiltonian invariants}
{\label{sct:the-inv}}

Before stating the main result of this section, we need 
some notations. We define the evaluation maps 
\begin{eqnarray}{\label{t-s-ev}}
\begin{array}{c}
\wtld\sev_k:{\mfrak X}_{B}\times {(P^c)}_o^k/
\crl{G}_0^c(P)\lar X^k,
\\[1ex]
\bigl[ (U,A),(p_1,\dots,p_k)\bigr] \lmt 
\bigl( U(p_1),\dots,U(p_k)\bigr)
\end{array}
\end{eqnarray} 
and 
\begin{eqnarray}{\label{s-ev}}
\begin{array}{c}
\sev_k:{\mfrak X}_B^{\rm s}\times {(P^c)}_o^k/
\crl{G}^c(P)\lar X^k,
\\[1ex]
\bigl[ (U,A),(p_1,\dots,p_k)\bigr] \lmt \bigl( 
U(p_1),\dots,U(p_k)\bigr)
\end{array}
\end{eqnarray}
and
\begin{eqnarray}{\label{bar-s-ev}}
\begin{array}{c}
\ovl\sev_k:{\mfrak X}_{B}\times {(P^c)}_o^k/
\crl{G}_0^c(P)\lar {({\cal P}_r^c(A_0)\times X)}^k,\\[.5ex] 
\bigl[ (U,A),\unl p\bigr] \lmt 
\bigl( [h(A-A_0),R_{f_A}\unl p],U(\unl p)\bigr).
\end{array}
\end{eqnarray} 
Computation \eqref{well-def} shows that the last 
evaluation is well-defined. All of them are $G^k$-equivariant 
and the last map induces on $M_{C,k}(\bar X;\bar B)$ the usual 
evaluation  
\begin{eqnarray}{\label{ev}}
\begin{array}{c}
ev_k:M_{C,k}(\bar X;\bar B)
={\mfrak X}_{B}\times {(P^c)}_o^k/\crl{G}_0^c(P)\times G^k
\lar {\bar X}^k.
\end{array}
\end{eqnarray}
The key for understanding the relationship between the 
analytic point of view and the algebraic one developed in 
the present paper is the diagram
\begin{align}{\label{diag-1}}
\xymatrix{%
{\mfrak X}_B^{\rm s}\times {(P^c)}^k_o/\crl{G}_0^c(P)
\ar[d]^-{\hbox{
\begin{minipage}[t]{20ex}
\scriptsize  
quot out the free\\ 
$G$-action (for $k\geq 1$)
\end{minipage}}}
&
\kern-1ex\subset\kern-5ex
&
{\mfrak X}_{B}\times {(P^c)}^k_o/\crl{G}_0^c(P)
\ar[d]^-{\hbox{
\begin{minipage}[t]{15ex}
\scriptsize  
quot out the free\\
$G^k$-action
\end{minipage}}}
\\
{\mfrak X}_B^{\rm s}\times {(P^c)}^k_o/\crl{G}^c(P)
\ar[d]^-{\hbox{
\begin{minipage}[t]{20ex}
\scriptsize 
quot out the (not necessarily free) 
$G^k$-action
\end{minipage}}}
&\kern-6ex &
{\mfrak X}_{B}\times {(P^c)}^k_o/\crl{G}_0^c(P)\times G^k
\ar@{=}[d]^-{\hbox{
\begin{minipage}[t]{15ex}
\scriptsize 
proposition \ref{prop:1-1}
\end{minipage}}}
\\ 
\kern2em 
S_{C,k}(X;B)\lrdash M_{C,k}(\bar X;\bar B)\invq G 
&
\ldashar\kern-5ex 
&
M_{C,k}(\bar X;\bar B).\kern2ex 
}
\end{align}
Since $G/K(={\mbb R}_{>0}^r)$ is contractible, $K$ and 
$G$-equivariant cohomologies of $X$ coincide; we shall 
prefer $G$-equivariant classes. Recall that $2D$ denotes the 
expected dimension of $S_{C,k}(X;B)$ and its formula is given 
by \eqref{eq:exp-dim}. 

\begin{definition}{\label{def:s-invar}} 
The Hamiltonian invariant introduced in \cite{cgs,mun} 
is defined in the following way: consider an equivariant 
coho\-mo\-lo\-gy class $\alpha\in {H_G^*(X)}^{\otimes k}$ 
with $\deg\alpha=2D$. Under the assumption that $G^k$ acts 
freely on ${\mfrak X}_B^{\rm s}\times {({\cal P}^c)}^k_o/
\crl{G}^c(P)$, the pull-back defines a cohomology class on 
$S_{C,k}(X;B)$ denoted the same. The invariant is 
\begin{eqnarray}{\label{eq:s-invar}}
\Phi^{C,k}_{X,B}(\alpha):=
\int _{S_{C,k}(X;B)}{(\sev _k)}^*\alpha.
\end{eqnarray} 
I have to say that $\Phi$ is defined this way only when 
$S_{C,k}(X;B)$ has the correct dimension. For this reason the 
authors in \cite{cgs,mun} work with perturbations of $S_{C,k}(X;B)$. 
In algebraic context, one should integrate over a $\pi$-relative 
virtual cycle, as mentioned in remark \ref{rk:probl}.
\end{definition}

It is conjectured in \cite{cgs} that for special choices of 
$\alpha$ and $B$ the invariant $\Phi$ should coincide with a 
Gromov-Witten invariant of $\hat X=X\invq G$. More precisely,
\vskip1ex 

\nit{\sl \unl{Con}j\unl{ecture}\,} {\it Take 
$\alpha\in {H_G^*(X)}^{\otimes k}$ and 
$\hat B\in H_2(\hat X;{\mbb Z})$, with $\hat X:=X\invq G$ 
and let $\hat\alpha\in H^*\bigl(\hat X^k\bigr)$ and 
$B\in H_2^G(X;{\mbb Z})$ be respectively the classes 
defined by 
$$
\hat X\srel{\sim}{\llar}
E\times_G X^\sst \lar E\times_G X.
$$ 
Then $\Phi _{C,k}^{X,B}(\alpha)
=GW_{C,k}^{\hat X,\hat B}(\hat\alpha)$.\vskip1ex}

Our goal is to prove this conjecture under the same transversality 
assumptions as in theorem \ref{thm:comparison}, when the invariant 
homology class $B\in H_2^G(X;{\mbb Z})$ is induced from $X^\sst$; 
the reason for this restriction was discussed in the end of the 
last section. So we will deal with topologically trivial 
$({\mbb C}^*)^r$-bundles over a smooth curve $C$. For 
${\cal P}_r^c\rar {(Pic^0C)}^r\times C$ the universal 
$({\mbb C}^*)^r$-bundle (trivialized at a point $\zeta_0\in C$), 
we define $\bar X:= {\cal P}_r^c\times_{({\mbb C}^*)^r}X$.

\begin{theorem}{\label{thm:compar2}}
Let the torus $G\cong ({\mbb C}^*)^r$ act on the irreducible 
projective variety $X$, and consider a linearization of this 
action in the very ample line bundle ${\cal O}_X(1)\rar X$. 
Denote $B\in H_2(X;{\mbb Z})$ a class which can be represented 
by a morphism $C\rar X^\sst$, where $C$ is a smooth projective 
curve with ${\rm Aut}(C)=\{{\rm id}_C\}$, and let 
$\hat B\in H_2(\hat X;{\mbb Z})$ be the class induced by 
the projection $\phi:X^\sst\rar X^\sst/G=\hat X$. Suppose that

{\rm (a1)\,} $G$ acts freely on the $G$-semi-stable locus of 
$X$, so that $X^\sst\rar\hat X$ is a principal $G$-bundle and 
denote $\bar B\in H_2(\bar X;{\mbb Z})$ the class induced 
as in lemma \ref{lm:u-bar};

{\rm (a2)\,} $\ovl M_{C,k}(\bar X;\bar B)$ is generically 
smooth and has the expected dimension;

{\rm (a3)\,} every irreducible component of 
$\ovl M_{C,k}(\bar X;\bar B)$ contains a point 
represented by a morphism 
$C\rar\bar X^\sst={\cal P}_r^c\times_{({\mbb C}^*)^r}X^\sst$;

{\rm (a4)\,} $M_{C,k}(\hat X;\hat B)\subset
\ovl M_{C,k}(\hat X;\hat B)$ is dense.

\nit Then for any $\alpha\in {H^*_G(X)}^{\otimes k}$, 
$$
GW^{C,k}_{\hat X,\hat B}(\hat\alpha)=\Phi^{C,k}_{X,B}(\alpha),
$$
where $\hat\alpha\in {H^*(\hat X)}^{\otimes k}$ is the class 
induced by $\alpha$.
\end{theorem} 

\begin{proof} 
Let me start explaining the guiding idea: we have learned in 
section \ref{sct:conseq1} that we should transfer the integrals 
on $\ovl M_{C,k}(\hat X;\hat B)$ used for defining Gromov-Witten 
invariants of $\hat X$ to integrals on 
$\ovl M_{C,k}(\bar X;\bar B)\invq G$ because 
$\bar X\invq G={(Pic^0C)}^r\times C\times \hat X$. 
All the evaluation maps involved for making these 
computations live on the right-hand-side of \eqref{diag-1}. 
On the other hand, the invariant $\Phi$ is defined using 
equivariant cohomology classes which are pulled-back to 
$S_{C,k}(X;B)$ by maps which live on the left-hand-side of 
\eqref{diag-1}. So, loosely speaking, what we have to do is 
to pass from the left to the right in \eqref{diag-1}. 

For making this passage, we need to understand the relationship 
between the various group actions and evaluation maps which 
appear in the context. On 
${\mfrak X}_B\times {(P^c)}^k_o/\crl{G}_0^c(P)$ there 
are two commuting actions. First, $G\subset\crl{G}^c(P)$ 
acts by constant complex gauges 
$$
g\times \bigl[ (U,A),\unl p\bigr] = 
\bigl[ (gU,A),R_g\unl p\bigr].
$$
The map $\wtld\sev_k$ is $G$-{\em invariant} for this action, 
while $\ovl\sev_k$ is $G$-{\em equivariant} for the diagonal 
right action of $G$ on the ${\cal P}_r^c$-factor of 
${({\cal P}_r^c\times X)}^k$.

Secondly, $G^k$ acts on 
${\mfrak X}_B\times {(P^c)}^k_o/\crl{G}_0^c(P)$ 
with quotient $M_{C,k}(\bar X;\bar B)$: 
$$
(g_1,\dots,g_k)\times
\bigl[ (U,A),(p_1,\dots,p_k)\bigr]
=\bigl[ (U,A),(R_{g_1}p_1,\dots,R_{g_k}p_k)\bigr].
$$ 
The evaluation map $\ovl\sev_k$ is $G^k$-equivariant for the 
$G$-action on ${\cal P}_r^c\times X$ on both terms.

\nit{\sl Convention\,} In what follows, the symbol ``$\sim$'' 
will denote homotopy equivalence and the letter ``$\jmath$'' 
obvious inclusions. For understanding better the forthcoming 
calculations, we should keep in mind that for integration 
purposes homotopy equivalent spaces are equal. 
\begin{align}{\label{eq:S}}
\begin{array}{cl}
S_{C,k}(X;B) &
\sim E^k\times_{G^k}\left( 
{\mfrak X}_B^{\rm s}\times {(P^c)}^k_o/\crl{G}^c(P)
\right)\\[1ex] 
&
\sim 
E^k\times _{G^k}\left( 
E\times_G\left( 
{\mfrak X}_B^{\rm s}\times {(P^c)}^k_o/\crl{G}_0^c(P)
\right)\right)\\[1ex] 
&
=
E\times_G\left( 
E^k\times_{G^k}\left( 
{\mfrak X}_B^{\rm s}\times {(P^c)}^k_o/\crl{G}_0^c(P)
\right)\right).
\end{array}
\end{align} 
The map 
$$
\ovl\sev_k:
{\mfrak X}_{B}\times {(P^c)}^k_o/\crl{G}_0^c(P)
\lar {({\cal P}_r^c\times X)}^k
$$ 
being $G^k$-equivariant and $G$-invariant, induces 
\begin{align}{\label{follow}}
\xymatrix{
E\times_G \left( 
E^k\times_{G^k}\left( 
{\mfrak X}_{B}\times {(P^c)}^k_o/\crl{G}_0^c(P)
\right)\right)
\ar[d]^-{\ovl\sev_k}\ar[r]^-{\sim}
&
E\times_G M_{C,k}(\bar X;\bar B)
\ar[d]^-{ev_k}
\\
E\times_G\left(
E^k\times_{G^k}{({\cal P}_r^c\times X)}^k
\right)
\ar[r]^-\sim
& 
E\times_G \bar X^k.
}
\end{align}
On the other hand, 
\begin{align}{\label{eq:X}}
\xymatrix{
E\times_G\left(
E^k\times_{G^k}{({\cal P}_r^c\times X)}^k\right)
\ar[d]^-{\epsilon}\ar[r]^-{\sim}
&
E\times_G \bar X^k
\\
E\times_G (E^k\times_{G^k} X^k)
=\BG\times(E^k\times_{G^k} X^k). 
& 
}
\end{align}
The reason for the last equality is that $G$ acts on 
${\cal P}_r^c$ only and therefore the $G$-action on 
$E^k\times_{G^k} X^k$ is trivial. The class $\alpha$ we start with 
lives in ${H_G^*(X)}^{\otimes k}$; pulling it back using $\epsilon$, 
we get the class $\bar\alpha\in H_G^*(\bar X^k)$. Relation 
\eqref{eq:S} implies that $\sev_k=\epsilon\circ\ovl\sev_k$ and we 
deduce from diagrams \eqref{follow} and \eqref{eq:X} that 
$\sev_k^*\alpha=ev_k^*\bar\alpha\in H_G^*(M_{C,k}(\bar X;\bar B))$. 

The invariant $\Phi$ can therefore be defined as 
\begin{align}{\label{transfer1}}
\Phi(\alpha)=\int_{\ovl M_{C,k}(\bar X;\bar B)\invq G}
\kern-2ex \jmath_{M^\sst}^*ev_k^*\bar\alpha,
\end{align} 
where the relevant maps fit in the diagram
\begin{align}{\label{diag-2}}
\xymatrix{
&E\times_G \ovl M_{C,k}(\bar X;\bar B)
&\kern-2.5em 
\raise1.5ex\hbox{$\xrar{\phantom{M}ev_k\phantom{M}}$}
\kern-2.5em 
& 
E\times_G \bar X^k
\\ 
&E\times_G\ovl M_{C,k}(\bar X;\bar B)^\sst
\ar[u]_-{\jmath_{M^\sst}}\ar[d]^-{\sim}
& &
E\times_G {(\bar X^k)}^\sst
\ar[u]_-{\jmath_{{(\bar X^k)}^\sst}}\ar[d]
\\
S_{C,k}(X;B)=\kern-3.55em
&\ovl M_{C,k}(\bar X;\bar B)\invq G
&\kern-2.5em 
\raise1.5ex\hbox{$\srel{\wht{ev}_k}{\longdashar}$}
\kern-2.5em 
&
\bar X^k\invq G.
}
\end{align}

From the diagram 
$$
\xymatrix{
 & & E^k\times_{G^k} X^k
\\
E\times_G 
\left(
E^k\times_{G^k}{({\cal P}_r^c\times X^\sst)}^k
\right)
\ar[d]^-{\sim}
\ar[rr]^-{\epsilon}
& & 
E^k\times_{G^k}{(X^\sst)}^k
\ar[u]_-{\jmath_{X^\sst}}
\ar[d]^-{\sim}
\\
E\times_G {(\bar X^\sst)}^k 
\ar[r]^-{\sim}
& 
{(\bar X^\sst)}^k/G
\ar[r]^-{q}
&
\hat X^k
}
$$
we see that $\jmath_{\bar X^\sst}^*\bar\alpha
=\jmath_{\bar X^\sst}^*\epsilon^*\jmath_{X^\sst}^*\alpha
=q^*\hat\alpha$, 
where $\hat\alpha$ is by definition the cohomology class on 
$\hat X^k$ determined by $\alpha$.

The assumption (a3) says that the subset 
$M_{C,k}(\bar X;\bar B)^o\subset
\ovl M_{C,k}(\bar X;\bar B)^\sst$ of morphisms with image contained 
in $\bar X^\sst$ is dense. We deduce, going the other way round in 
\eqref{diag-2}, that 
\begin{align}{\label{eq:final}}
\disp\Phi(\alpha)=\int_{\ovl M_{C,k}(\bar X;\bar B)\invq G}
\kern-2ex{(\wht{ev}_k)}^*q^*\hat\alpha.
\end{align}
At this point we have finally moved from the left-hand-side to 
the right-hand-side of \eqref{diag-1} as we wished. Notice that 
since the maps involved in the computations are rational, the 
pull-backs are defined as in \ref{def:corresp}. 

The composition of morphisms $\bar u:C\rar \bar X$ representing 
the class $\bar B$ with the projection $\bar X\rar{(Pic^0C)}^r
\times C$ is of the form $\{\tau\}\times{\rm id}_C$ because 
$\bar B$ induces the class $0\oplus [C]$. Consequently, a map 
$\bar u:C\rar \bar X^\sst$ defines a map $C\rar C\times\hat X$ 
which is the identity on the first component.

\begin{lemma}{\label{lm:birtl}} 
The map 
$$
T:\ovl M_{C,k}(\bar X;\bar B)\invq G\dashar
\ovl M_{C,k}(C\times\hat X;[C]+\hat B)
$$ 
is birational. 
\end{lemma}

\begin{proof}
Assumption (a4) implies that $T$ is dominant if every morphism 
$\hat u:C\rar\hat X$ representing $\hat B$ is in its image. From 
the diagram 
$$
\xymatrix{
\hat u^*X^\sst
\ar[d]\ar[r] &
C\times X^\sst 
\ar[d]\ar[r] &
X^\sst
\ar[d]^-{\phi}
\\ 
C\ar[r]^-{{\rm id}_C\times\hat u}&
C\times \hat X
\ar[r] &
\hat X
}
$$
we deduce that the pull-back $\hat u^*X^\sst\rar C$ is a 
topologically trivial, holomorphic principal bundle (this 
follows from the assumption that $\hat B$ is induced by a 
class $B\in H_2(X^\sst;{\mbb Z})$); it is therefore isomorphic 
to ${\cal P}_\tau:={{\cal P}_r^c}|_{\{\tau\}\times C}$ for a 
certain $\tau\in{(Pic^0C)}^r$. This data induces the 
map $\bar u:C\rar {\cal P}_\tau\times_G X^\sst
\subset {\cal P}_r^c\times_G X=\bar X$ which represents 
the class $\bar B$ and also $T(\bar u)=\hat u$.  

Now let us prove that $T$ is generically injective. Using 
assumption (a3), we may restrict ourselves to the (dense) 
open subset $M_{C,k}(\bar X;\bar B)^o$ representing morphisms 
whose image is contained in $\bar X^\sst$. Let us assume that 
$$
\begin{array}{ccccc}
C& \srel{\bar u}{\lar}& 
{\cal P}_\tau\times_G X^\sst& & \\[-1ex]
 & & &\kern-2ex \sear^\phi& \\
 & & & &\kern-4ex C\times\hat X\\
 & & &\kern-2ex \near_\phi& \\[-1ex]
C& \srel{\bar u'}{\lar}& 
{\cal P}_{\tau'}\times_G X^\sst& & 
\end{array}
$$
are such that $\phi\circ \bar u=\phi\circ\bar u'$. Then 
${(\phi\circ \bar u)}^*(X^\sst\rar\hat X)=
{(\phi\circ \bar u')}^*(X^\sst\rar\hat X)$ and lemma 
\ref{lm:tau} implies that $\tau=\tau '$. Then 
$C\xrar{\bar u,\bar u'}{\cal P}_\tau\times_G X^\sst$ 
induce the same map to $\hat X$ and consequently for any 
$\zeta\in C$, $\bar u'(\zeta)=g(\zeta) u(\zeta)$ for a unique 
$g(\zeta)\in G$. The morphism $C\rar G$ must be constant, so 
$\bar u'=g\bar u$ and they define the same point in 
$\ovl M_{C,k}(\bar X;\bar B)\invq G$.
\end{proof}

Since $C$ has trivial automorphism group, 
$$
\ovl M_{C,k}(C\times\hat X;[C]+\hat B)
\srel{1:1}{\longlrdashar}
\ovl M_{C,k}(\hat X;\hat B).
$$

Finally, from the diagram  
$$
\begin{array}{ccc}
\ovl M_{C,k}(\bar X;\bar B)\invq G\kern-8ex&
\srel{\wht{ev}_k}{\hdashpiece\kern-2.5pt\dashar} 
{\bar X}^k\invq G
\srel{q}{\hdashpiece\kern-2.5pt\dashar}&
\kern-1.5ex\hat X^k\\ 
_{1:1}\vdashto\kern-4ex& &\kern-2.5ex\vequal{3mm}\\ 
\ovl M_{C,k}(C\times\hat X;[C]+\hat B)\kern-1.7ex&
\xrar{{ev}_k^{\hat X}}
{(C\times\hat X)}^k
\xrar{{\rm pr}_{\hat X}}&
\kern-1.5ex\hat X^k
\end{array}
$$
we conclude that the invariant $\Phi$ does coincide with 
a Gromov-Witten invariant of $X$. Indeed, 
$$
\begin{array}{cl}\disp
\Phi^{C,k}_{X,B}(\alpha)\kern-1ex&\disp
=\int _{\ovl M_{C,k}(\bar X;\bar B)\invq G}
\kern-2ex\wht{ev}_k^*q^*\hat\alpha
=\int _{\ovl M_{C,k}(C\times\hat X;[C]+\hat B)}
\kern-2ex{({ev}_k^{\hat X})}^*{\rm pr}_{\hat X}^*\hat\alpha
\\[1.5ex]
&\disp
=\int_{\ovl M_{C,k}(\hat X;\hat B)}{({ev}_k^{\hat X})}^*\hat\alpha.
=GW_{\hat X,\hat B}^{C,k}(\hat\alpha).
\end{array}
$$
\end{proof}

\nit I would like to conclude with the remark:

{\sl Why should be interesting the Hamiltonian invariants?\,} 
I have mentioned in the introduction that the starting point of 
this study was the problem of comparing the GW-invariants of a 
quotient with the invariants of the variety we start with. Simple 
dimensional counting shows that --except in genus zero--  the 
question is not well-posed in this form: the dimension of the space 
of morphisms from curves to a quotient variety is larger than the 
dimension of the space of morphisms into the starting variety. The 
difference between these dimensions is exactly the dimension of the 
moduli space of principal $G$-bundles over a curve; this can be 
explained noticing that, for morphisms $v:C\rar\hat X$ which 
represent an {\it a priori} given homology class, the holomorphic 
type of the the pull-backs $v^*X^\sst\rar C$ changes within a fixed 
topological type. In this way, the space of principal bundles 
over curves with fixed topological type naturally enters into the 
scene. 

Equality \eqref{eq:final} brings our attention to another aspect 
of the problem: GW-invariants of $\hat X$ can be computed (under 
suitable transversality conditions) in the following way
$$
GW_{\hat X,\hat B}^{C,k}(\hat\alpha)
=\int_{\ovl M_{C,k}(\bar X;\bar B)\invq G}\kern-1em
{(q\circ\wht{ev}_k)}^*\hat\alpha,
$$
for $q:{(\bar X^\sst)}^k/G\rar \hat X^k$. The interesting part is 
that there is a natural projection 
$$
\pi:\ovl M_{C,k}(\bar X;\bar B)\invq G\rar {(Pic^0C)}^r,
$$
so that one can further write 
$$
GW_{\hat X,\hat B}^{C,k}(\hat\alpha)
=\int_{{(Pic^0C)}^r}\kern-1em 
\pi_*{(q\circ\wht{ev}_k)}^*\hat\alpha
$$
and now the integration takes place on the Picard torus of the 
curve. One may hope that enumerative invariants of suitably chosen 
$\hat X$'s can be expressed in terms of intersection numbers on 
the Picard variety.


\end{document}